# Galois Got his Gun

*Frédéric Brechenmacher*

Univ. Lille Nord de France, F-59000 Lille.

U. Artois, Laboratoire de mathématiques de Lens (EA 2462),

rue Jean Souvraz S.P. 18, F- 62300 Lens, France

frederic.brechenmacher@euler.univ-artois.fr

> Tout cela étonnera fort les gens du monde, qui, en général, ont pris le mot Mathématique pour synonyme de régulier. Toutefois, là comme ailleurs, la science est l'œuvre de l'esprit humain, qui est plutôt destiné à étudier qu'à connaître, à chercher qu'à trouver la vérité […]. En vain les analystes voudraient-ils se le dissimuler : quand ils arrivent à la vérité, c'est en heurtant de côté et d'autres qu'ils y sont tombés.
>
> Évariste Galois, « Discussion sur les progrès de l'analyse pure» [Galois 1962, p. 15]. Quoted in [Boutroux 1920, p. 190], [Picard 1924a, p. 31], [Picard 1924b, p. 327], [Picard 1925, p. 23].

## Introduction

When French students in mathematics were called to arms in World War I, the long dead Évariste Galois also got his gun. The latter was indeed both a heroic model for the ghosts of the *École normale supérieure* who had perished in the trenches before they had a chance to prove their quality [Julia 1919],[1] and an inspirational scientific figure for the mathematicians involved in the *Commission de balistique* [Hadamard 1920, p. 438]. After 1916, Émile Picard, one of the prominent authorities of the *Académie des sciences*, celebrated Galois as one the heroes of the universality of the French style of thinking, i.e., as one of the *Grands savants* whose accomplishments were raised as shields against the manifest of the ninety-three German scientists.

These memorial references to Galois during World War I echoed the episode when Galois himself took his gun in 1832. The episode of the duel is indeed a highly inspirational symbol of the tension between the realm of pure mathematical ideas and the physical world with its dramatic, political, and, actually, quite plain realities. This episode both escapes the boundaries of generic categories and, at the same time, delimits such boundaries. As a matter of fact, one of the main specificities of the public references to the figure of Galois in the long run is, one the one hand, their great generic diversity and, on the other hand, the boundary-works these references involve across fiction and history, (pure) mathematics and politics, biography and literature, academic authority and public opinion, etc. [Albrecht and Weber 2011].

---

[1] Quoted and transalted in [Aubin forhcoming].

The present paper, in a way, appeals to the figure of Galois for investigating the gates between mathematics and their "publics." Whether it is opened or closed, a gate is always delimiting the inside from the outside. The very large "public" dimension of Galois as a "mathematician" thus places the latter in the position of a gatekeeper. The figure of Galois draws some lines of/within mathematics for/from the outside and these lines in turn sketch the silhouette of Galois as a historical figure. There is indeed a striking contrast between the long-term simplicity of the main categories involved in most public discourses on Galois's achievements (e.g., "equations," "groups," "algebra," "France," "Germany") and the variety of actual references to Galois in mathematical papers. More precisely, the recurrent public statements that Galois's works were mainly concerned with the groups of general equations involved in the "Mémoire sur les conditions de résolubilité des équations par radicaux" (the *Mémoire* for short) contrast with the new light some recent researches in the history of mathematics have shed on the actual reception of Galois's works in the 19$^{th}$ century.[2]

The main categories involved in the usual presentations of Galois's mathematical achievements cannot be disconnected from their public dimensions. Moreover, despite their apparent long-term stability, these categories have actually taken on various meanings. Let us consider a few examples. In 1962, opening Robert Bourgne and Jean-Pierre Azra's critical extended edition of Galois's collected works, Jean-Dieudonné made the statement that it would be superfluous to "say once again after so many others what mathematics owes to Galois," because "everyone knows that his (Galois) ideas are the source of modern Algebra itself" [Galois 1962, p. 1]. This public discourse has underlying it the role of authority Dieudonné endorsed in expressing the impressive ("the ideas at the source of") self-evidence ("everyone knows") of the relevant categories for identifying Galois's works ("modern Algebra").

But let us now compare Dieudonné's introduction to that of Picard in the 1897 reprinting of Galois's works. Picard also celebrated Galois for his introduction of the notion of group. But he nevertheless opposed the figure of Galois to modern algebra in insisting on the analysis of groups of operations. As a matter of fact, some authoritative figures in French mathematics successively advanced two claims in regard with Galois, equations, and groups: one for the French style of thinking in analysis as opposed to the German arithmetic and algebra (1900-1930); and a second claim, symmetrically, that celebrated German conceptual algebra as opposed to the older French computational approaches (1930-1970).

But the public dimension of Galois is not limited to academic celebrations. Galois was indeed already a public figure as a mathematician early on in the 1830s, i.e., much before Liouville rehabilitated his mathematical works in 1846. Moreover, Galois's public figure shows some long-term continuity in parallel – and, actually, in spite of – the discontinuous interventions of mathematical authorities.

---

[2] See [Ehrhardt 2007] for a detailed analysis of the difficulties raised by the traditional history of the reception of Galois's works in regard with the development of group theory and Galois theory.
The actual reception of Galois's works especially involved the special equations associated to elliptic functions (modular equations) [Goldstein 2011], number-theoretic imaginaries in connection to the analytic representation of substitutions [Brechenmacher 2011], and some analogies between algebraic and differential equations [Archibald 2011].

These continuities between public and academic discourses raise the difficult methodological issue of the circulations, or the transfers, between the various fields of knowledge or activities that are involved in some highly heterogeneous corpora. These often include altogether mathematical papers, novels, newspapers, dictionaries, academic publications, etc. This issue has been usually investigated in connection to the one of the transfers of symbolic capital (e.g. from sciences to politics, from mathematics to philosophy etc.). But in the present paper I would rather aim at discussing the circulations of the collective categories that are used in discourses on mathematics. It is thus my aim to investigate how the lines of delimitations of some fields may themselves be transferred from a field to another, such as for instance from mathematics to the public sphere or from the public sphere to the historiography.

Let us consider the example of Picard who claimed in 1883 he aimed at approaching differential equations in analogy with Galois theory of algebraic equations. This claims marks the starting point of differential Galois theory. But it is striking that in the 1870s, several authors had already developed various analogies between differential and algebraic equations without referring to Galois at all.[3] It is thus not impossible that Picard was at least partially aiming at taking a public stand in reacting to the increasing uses of Galois's name by Klein and his followers in Germany. Moreover, it is striking that in the 1890s, Picard, Vessiot and Drach appealed to some "analogies" to Galois's works very different from one another [Archibald 2011]. Was there any echo of the inspirational power of the public figure of Galois in these specialized mathematical works? Shall one see any connection between the diversity of the analogies involved in differential Galois theory and the variety of "generalizations" of Galois theory that were developed by some French mathematicians at the time of World War II?[4]

Caroline Ehrhardt has given a detailed analysis of how Galois was made an icon of mathematics. The present paper lays the emphasis on the different issue of the categories that have been used in various types of discourses on Galois's works. I shall therefore follow some of Herbert Mehrtens's perspectives on public discourses on mathematics in Germany in the first half of the 20th century. Mehrtens has indeed shown that even though public discourses did usually not directly reflect contemporary mathematics, they were nevertheless in a complex interaction one another, as the two sides of a same coin [Mehrtens 1990]. I will nevertheless not appeal to Mehrtens's antagonistic dualism between the modernism of abstract mathematics, as heralded by actors such as David Hilbert, and the counter modernism focus on intuition and on the applications of mathematics, as heralded by mathematicians such a Felix Klein.

Mehrtens's dichotomy has underlying it the thesis that modernism is related to internationalism and to the liberal attitude of autonomous professionals, while counter modernism lends itself to non-scientific ideological integration of rather disparate values and identification, and is thus correlated to nationalism and eventually also to racism. But this dualism does not reflect the variety of the roles played by actors such as Picard who were promoting both the essentialist view that there is a natural substance to the truth and meanings of mathematics (counter-modern) and the creativity of mathematicians or the importance of pure mathematics in regard with applications (modern). Some peripheral actors such as Jean-Armand de Séguier also show the limits of Mehrten's duality: while Séguier's works was dominated by abstraction (modern), the latter connected the beauty of pure abstraction to the "sublime" in a religious sense (counter- modern; recall that Séguier was a Jesuit abbot).

---

[3] Among these authors, Jordan can hardly be suspected of not having studied Galois's works closely enough.

[4] "Generalizing" Galois theory was one of the main topics of the papers published by the members of the Bourbaki group at the time of World War II. The meanings of such generalizations varied from the consideration of non-commutative fields to the one of infinite extensions. When the Dubreil-Châtelet seminar on algebra and number theory was created at the *Faculté des sciences de Paris* in 1947, six of the seven talks were connected to "generalizations" of "Galois theory".

In the present paper, I shall especially focus on the period from 1890 to 1930 in France. This time-period opens with Galois's entrance in the mathematical pantheon and closes when discourses on the universality of the French style of thinking were at a climax. In the meantime, some French authorities repeatedly opposed the figure of Galois to the idea that algebra and arithmetic should have some autonomous value in regard to analysis. That Galois's works have been publicly inscribed in the analysis for decades has fallen into oblivion after the development of algebra as both an autonomous discipline and a professional specialty in the 1930s-1950s.

Because of the long-term public dimension of the figure of Galois, we shall also appeal to a multi-scale analysis in considering from the perspective of the *Belle Époque*'s Galois a few key episodes in the broader time period 1830-1940.

As an actor's category, the notion of "public" was actually much involved in discourses on Galois since the 1830s. It especially played a key role in the 1890s when Galois was celebrated for having introduced a general notion of group in the special case of algebraic equations. Most contemporary discourses on Galois indeed revolved on the tension between Galois individual "ideas" and their collective dimensions. More precisely, the presentation of Galois's "influence on the developments of mathematics" [Lie 1895] was structured by the tension between the "obscurity" of Galois's original ideas and the illumination provided by Camille Jordan's unfolding of the group-theoretical nature of Galois's ideas,[5] thereby making these ideas "intelligible to the general public." [Pierpont 1897, p.340]. In the 1890s, most histories of the theory of equations adopted a three-act structure: before Galois, Galois, and how Jordan's 1870 *Traité des substitutions et des equations algébriques* had "made a knowledge of Galois's theory possible to all the world" [Pierpont 1897, p.340].

To be sure, the epithet public took on various meanings in connection with mathematics in various time-period and socio-political contexts. We shall thus certainly not appeal to a reification of the notion of public to distinguish between mathematical and public discourses. Here, a very simple criterion actually allows a dynamical distinction between the changing categories of public and mathematics: in the present paper, we shall designate as "public discourses on Galois" all the texts that are appealing to an authority such as Liouville,

---

[5] [Ehrhardt 2007, p. 1-45] has historicised the category of the "intelligibility" of Galois's writings.

Picard, or Dieudonné.[6] Other texts referring to Galois will be designated as mathematical papers.

This criterion is very efficient. Actually, authorities are often explicitly quoted word-for-word in public discourses while they are not even alluded to in mathematical papers. It may be added that a second characteristic of general public discourses is that no specific reader or listener is known. In regard to our two characteristics, some textbooks are clearly on the borderline between public discourses and mathematical papers. Serret's *Cours d'algèbre supérieure* is a typical example of a book that quoted Liouville word-for-word on the one hand, and of which some individual readers are known on the other hand.

A tension between individual and collective dimensions of mathematics thus lies at the roots of the distinction between public discourses and mathematical papers. This tension provides the opportunity to investigate some mathematical *personae,* in the sense of the socio-cultural identities expressed by the various roles that have been taken on by some individuals (such as Picard) and that have been assigned to some others (such as Jordan).[7] These roles especially raise the issue of the relational nature of certain discourses on mathematics. To be sure, eulogies both aggrandize the dead and crown the living beings.[8] But the issue of the relational nature of discourses on Galois is not limited to the reflective nature of eulogies. Indeed, very few historical sources are available on Galois's life. By filling the holes in Galois's biography, history has often met fiction and, conversely, fiction has often claimed its historical relevance [Albrecht & Weber 2011]. Some actors have thus been portraying themselves as well as the collective roles they were playing through the relations they have established to Galois.[9]

This issue of the self-portrait is much connected to the main actors' conscious perception that the parallel growth and fragmentation of both the public and the practitioners of mathematics made it compulsory for papers to target specific audiences and for authorities to express their presence publicly.[10] A paper published in *Le Figaro* in 1912 exemplifies this situation. The author, Louis Chevreuse, depicted with irony a world in which the status of manuscripts had turned from a promise of potential masterpieces to the "mysterious horror" that plagued contemporary play directors, novelists, editors, scientists etc. Cheuvreuse appealed to the "high authority" of Darboux who, as one of the two *secrétaires*

---

[6] The role played by authority in sciences has been emphasized in [Elias 1982]. Official discourses of authority at the turn of the century have been recently analyzed in [Weber 2012]. On the institutional authorities of French mathematics, see [Gispert 1991 & 1995] and [Zerner 1991] for the period from 1870 to 1914 and [Gispert & Leloup 2009] for the period from 1918 to 1939.

[7] On the various models of mathematical lives during the interwar period, see [Goldstein 2009]. On the interconnections of the notions of roles and persona, see [Aubin & Bigg, 2007]. On the notion of persona itself, see [Daston & Sibum, 2005], [Mauss, 1938].

[8] Cf. [Bonnet 1986] and [Weber 2012]. See especially de Parville's comment on Berthelot's first eulogy (of Lavoisier) as one of the two *secrétaires perpetuels* of the *Académie* in 1899: "now we know well Lavoisier and we known even better M. Berthelot" (quoted from [Weber 2012]).

[9] See [Albrecht & Weber 2011]: the figure of Galois calls for some analysis on the relational nature of biographic narratives. On the hermeneutical approach to biography and on the resulting instability of historical lives, see [Rupke 2008].

[10] This situation has been recently analyzed by Laura Turner in the case of the various roles taken on by Mittag-Leffler. See especially Turner's analysis of the way Mittag-Leffler contrasted the international public of *Acta Mathematica* with the local public of Scandinavian scientific congresses. See also the ways both Peano and Mittag-Leffler commented on the creation of an international journal as the simultaneous creation of a public [Turner 2011].

*perpetuels*, was in charge of opening the sealed envelopes of the manuscripts sent to the *Académie*. "But the [academic] sessions are short and the memoirs are long. And for all his mathematical genius, M. Darboux cannot make the part larger than the whole." [Cheuvreuse 1912, p. 1]. Darboux, the chronicler reported, thus sadly sighted at the view of the hundreds of unopened letters that encumbered the *Académie*'s archives. Picking an enveloppe randomly in a somewhat vain an desperate gesture, the academician discovered a memoir submitted forty-three years earlier by a mathematician named Gaston Darboux… Cheuvreuse then recalled the stories of Abel and Galois and concluded that "the *savants* who do not have a mind for getting old are doomed to be proven wrong." At first sight Cheuvreuse's claim might seem absurd in regard with Abel's and Galois's fames. But it actually ironically laid the emphasis on the relational nature of the contemporary discourses on Galois, whose authors were much more expressing themselves than presenting Galois's original ideas.

Some self-portraits have indeed been sketched as a result of the boundary-works that most references to Galois involve: academicians versus teachers, researchers versus amateurs, history versus fiction, mathematics versus politics, algebra versus analysis, France versus Germany, etc. As shall be seen in this paper, the categories involved in discourses on Galois shed light on some changing models of mathematical lives. These categories also highlight that the public expressions of such models have underlying them a long-term tension between academic and public legitimacies.

In sum, this paper will neither focus on the episode when Galois took his gun nor on later echoes of this episode. On the contrary, it aims at shedding light on the boundaries some individuals drew by getting Galois his gun. Robert d'Adhémar's obituary of Jordan exemplifies the kind of issues I would like to investigate here. This obituary indeed highlights how the intellectual filiation to Galois aggrandizes Jordan from a man to a great mathematician, thereby simultaneously setting a hierarchy of disciplines in a tension between the glorious world of "mathematical beings" and the plain sadness of the real world:[11]

> In 1860, Jordan was already devoting himself to the Algebra of order, i.e., an Algebra of ideas which is much higher than the Algebra of computations. He naturally followed Galois's works, this genius and disappointing child who was wounded in a ridiculous duel in 1832 and who consequently died at the age of 21. […] Jordan's discoveries were published in 1870 in the *Traité des substitutions et des equations algébriques*, which, after Abel and Galois, marks an immense progress in Algebra. […] Whenever Jordan manipulates a mathematical being, it is with the austere hold of his powerful claw. Wherever [Jordan] passes, the trench is cleared. [Adhémar 1922, p. 65].[12]

---

[11] Recall that Jordan had lost three sons and one grandson during the war, all killed in action. On the brutalization of the language in mathematics in connection to the violence experienced by the younger generation of mathematicians involved in WWI, see [Aubin forthcoming].

[12] Jordan s'applique, dès 1860, à l'Algèbre de l'ordre, l'Algèbre des idées, bien plus haute que l'Algèbre des calculs, et, tout naturellement, il continue l'oeuvre de cet enfant genial et décevant, Galois, qui, blessé dans un duel ridicule, mourut en 1832, âgé de 21 ans. […] Ses découvertes ont été publiées en 1870 dans le *Traité des substitutions et des equations* algébriques, qui marque, après Abel et Galois, un progress immense de l'Algèbre. […] Chaque fois qu'il manie un être mathématique, Jordan met sur lui sa griffe puissante et austere. Là où il a été, la tranchée est nettoyée !

# 1. From equations to groups: Galois goes public

The focus on groups and equations in the usual history of Galois's posterity is partly a consequence of the evolutions of mathematics, and especially of the increasing role played by object-oriented disciplines in the organization of mathematics (e.g. group theory, Galois theory, algebraic number theory, linear algebra etc.). It is well known that the writing of the history of a field characterizes the existence of such a field [Bourdieu 1976, p. 117]. This phenomenon makes it is difficult to access the collective organizations of knowledge that existed prior to the ones that are contemporary to us.[13]

Because of the retrospective concerns for the origins and diffusions of the key notions of nowadays Galois's theory, the roles played by some actors such as Charles Hermite and Leopold Kronecker have been recurrently underestimated. These works indeed hardly fit the object-oriented focus of nowadays disciplines [Goldstein 2011]. In parallel, Jordan's works have been celebrated for their role in the empowerment of group theory as a discipline in its own right. These works have thus been connected to those of Galois in an exclusive relation. As a result, Jordan's *Traité* has been cut from the collective dimensions in which it originally made sense [Brechenmacher 2011].

But the evolutions of public discourses on Galois are not limited to the impact of the evolutions of mathematics on the history of mathematics. Such discourses also depended on the agendas of their authors in targeting some audiences. To be sure, these agendas were not limited to the aim of popularizing contemporary mathematics. Aggrandizing men was one of the traditional roles of official eulogies. We shall thus start our investigations by looking into how the celebrations of the Great mathematician Galois at the turn of the 20th century interlaced mathematical categories, such as groups and equations, and some collective epistemic and moral values.

## 1.1. The public dimension of the passage from equations to groups

The rhythm of the official celebrations of Galois is chanted by the introduction Liouville, Picard, and Dieudonné successively gave to three editions of Galois' works (1846; 1897; 1962). Let us take an overview of the impulse given by these editions on the categories used in public discourses from 1830 to 1960.

In his *Avertissement* to the 1846 edition of Galois's works, Liouville rehabilitated Galois's *Mémoire* with regard to the doubts Lacroix and Poisson had expressed when they had reported on the memoir in 1831 [Ehrhardt 2010].[14] In doing so, Liouville claimed that Galois had laid the grounds for a "general" theory of the solvability of equations by radicals. While he did not specify the content of such a "general theory," Liouville celebrated the generality of "Galois's method" through its particular use in the proof of the *Mémoire*'s concluding theorem, i.e., the criterion that an irreducible equation of prime degree is solvable by radicals if and only if all of its roots are given by a rational function of two of them (the *Galois criterion* for short).

---

[13] Several recent case studies in the history of mathematics have shed light on various aspects of such issues. See [Goldstein 1999], [Goldstein and Schappacher 2007a], [Brechenmacher 2007], [Gauthier 2009], [Brechenmacher 2011], [Goldstein 2011].

[14] As is well known, Galois presented three different papers on the solvability of algebraic equations to the *Académie* between 1829 and 1831. The first two were lost and the third was returned with a request for clarifications [Ehrhardt, 2010].

The presentation of the *Galois criterion* as a particular application of a general theory of equations dominated public discourses on Galois's works until the mid-1890s. Liouville's presentation of Galois was in fact reproduced word for word in publications targeting larger audiences than specialized mathematical journals, e.g., the 1848 biography of Galois in the *Magasin encyclopédique* and the many notices that were published afterward in encyclopedic dictionaries. Moreover, the citation of Liouville citing Galois was also to be found in Serret's *Cours d'algèbre supérieure.* Despite the fact that the first edition of 1849 had made almost no use of Galois's works, its introduction nevertheless presented the *Galois criterion* as the endpoint of a *longue durée* history of the "theory of equations" involving Cardano, Lagrange, Ruffini, and Abel among others [Serret 1849, p. 1-4].

In the 1860s, Joseph Bertrand, then at the peak of his power on both the Academy and the Parisian mathematical world [Zerner 1991], also gave a public presentation of Galois's works. This presentation was different from the one of Liouville's *Avertissement* and was connected to the new presentation of "Galois theory of general equations" in the 1866 edition of Serret's *Cours*. We shall return to Bertrand later. For now, let us emphasize the fact that Liouville's 1846 presentation nevertheless continued to dominate public discourses on Galois's works until the 1890s, at which point it suddenly disappeared. Indeed, the main difference between the 1897 reprinting of Galois's works and the 1846 edition was that Picard's introduction was substituted to Liouville's *Avertissement*.

The 1897 reprinting was published two years after the celebrations of the centenary of the *École normale supérieure* when, one the one hand, Sophus Lie had lectured on "Galois's influence on mathematics," and, on the other hand, the historian Paul Dupuy had published a biography of Galois [Ehrhardt 2007, p.628-649]. While he did not even make the slightest allusion to the *Galois criterion,* Picard followed the role Lie had assigned to Jordan as the one who had "clarified, developed, and applied" substitution groups in regard to Galois's works on the solvability of equations [Lie 1895, p.4].[15] Jordan was thus presented as the immediate follower, the one who was in direct contact with Galois' ideas and who had generalized the latter's distinction between simple and compound groups to the notion of composition series [Picard *in* Galois 1897, p.viii].

Jordan was thus assigned a major role in the three-act story of the "predecessors", the "origins", and the "influence" of "Galois' ideas" [van der Waerden 1985, p.76-133]. He was presented as the one who mediated some individual ideas and the collective appropriation of a comprehensive theory. As has already been mentioned before, the role of go-between assigned to Jordan highlights that publicity was an important issue in the attribution to Galois of the passage from equations to groups. All reviews indeed emphasized that neither the reprinting of Galois's works nor Dupuy's biography were for the eyes of mathematicians only, including the reviews of journals with editorial orientations as different as the Jesuit *Études* and the Republican *La révolution française. Revue d'histoire moderne*. In the *Bulletin*

---

[15] Before Picard and Lie, Klein had already played an important role in the presentation of Galois as one of the founders of group theory.

*of the American Mathematical Society,* Yale professor James Pierpont regretted the absence of historical comments that would have helped the non specialist "reader who is becoming interested in Galois's theories" by pointing to "some of the principal differences between Galois's exposition of his theory and that received today." A reprinting of Dupuy's biography, Pierpont added, could have helped to make it "accessible to everyone."

### 1.2. On Picard's public authority

Picard was to become one of the most prominent authorities, not only of French mathematics but also of French science. In the 1920s, Picard had become a "distant divinity",[16] even though he was still much involved in the organization of both national and international science. He was also one of the public official representatives of science in the medias - a situation one may compare to others such as the one of Mittag-Leffler who was regularly featured in Swedish newspapers as an illustrious sage, not only of science but also of the unity of the Scandinavian countries [Turner 2011]. Let us reproduce the presentation by which *Le Figaro* introduced a short paper by Picard on the 3rd of May 1926:[17]

> Both a mathematician and a philosopher, Emile Picard is professor at the Sorbonne, elected at the *Académie française* and one of the *secrétaire perpétuel* of the *Académie des sciences*. His universal reputation results from his works in mathematical analysis which has drawn to Picard students from the entire world. Thanks to analysis, Picard achieved a supreme mastery on mathematics. It was thus from a higher perspective that he later entered into other sciences. In a series of usually concise writings, Picard has expressed such a great number of original thoughts and deep observations on the various fields of human activities - especially on physics - that a methodical scholar could find here materials for a whole new philosophy of sciences.[18]

One may add that in the 1920s Picard was also the head of the International Research Council and that he was also involved in the direction of most important French mathematical journals. Moreover, when they alluded to Picard's scientific achievements, the medias always referred to the latter's generalisation to differential equations of "Galois's fantastic discoveries" [Adhémar 1924]. We shall thus investigate further the connections between the public figure of Picard and the role the latter took on in the publicity of Galois's "ideas".

As shall be seen in greater details later, the public dimension of the figure of Galois had been challenging the authority of the *Académie* since the 1830s. In this context, both Liouville's 1846 edition and Picard's 1897 reprinting are landmarks of the Academy's progressive loss of monopoly on discourses on sciences. Both indeed highlight the increasing autonomy of mathematics as well as the increasing role played by professorship as a model of mathematical life. The 1846 edition in the *Journal de mathématiques pures et appliquées* had been immediately celebrated in the *Nouvelles annales de mathématiques*, a journal whose contributors were mostly teachers and their pupils. Moreover Olry Terquem, one of the editors of this journal, almost immediately opposed Galois's genius to the second-rate mathematical quality of the *École polytechnique*'s examiners, despite the central role this institution played in the mathematical world at the time [Terquem 1849]. As for the 1897 reprinting, it was supported by the *Société mathématique de France* while the celebrations at the *É.N.S.* had consecrated the leading role the professors who had been trained in this institution played in the elite of mathematics in France.

---

[16] According to Mandelbrojt, quoted in [Gispert and Leloup 2009, p. 46].

[17] As other dailies of the *bourgeoisie*, *Le Figaro* often reported on scientific issues on the model of the *rubriques mondaines*. Here, Picard's paper was concerned with the relations between science and religion (Picard himself was a devout catholic).

[18] Mathématicien et philosophe, professeur à la Sorbonne, membre de l'Académie française et secrétaire perpétuel de l'Académie des sciences. L'analyse mathématique a fait la réputation universelle de M. Émile Picard, et lui a attiré des élèves de tous les pays du monde. Parvenu, grâce à elle, à une maîtrise suprême dans les mathématiques, il a ensuite abordé les autres sciences, mais par le haut, si l'on peut dire. Dans un ensemble d'écrits généralement courts, il a exprimé sur les divers champs d'activité de l'esprit humain et particulièrement sur la physique, tant de pensées originales et de si fines observations, qu'un érudit méthodique y puiserait facilement les matériaux de toute une nouvelle philosophie des sciences.

But the autonomy of the community of professional mathematicians was nevertheless still a relative one. Recall that Liouville and Picard were both academicians. Both appealed to two of the academy's symbolic powers that had been built again the tribunal of public opinions, i.e., the power of authorizing a form of language [Bourdieu 2001] on the one had, and the power of deciding of the chosen ones - those who deserve to be aggrandized, on the other. The traditional rhetoric of the academic eulogies was indeed extending its shadow on both Picard's and Lie's discourses.

Lie, especially, had described Galois as "*cet immortel normalien.*" In 1911, *La revue de Paris* published an obituary of Jules Tannery that illustrates further the transposition to the *É.N.S* of the epistemic continuity modeled on the academician's immortality. Recall that Tannery also had played a role in celebrating Galois at the turn of the century [Tannery 1909].[19] According to *La revue de Paris*, Tannery's "two fold involvement (or duty) in truth," i.e., in both science and the teaching of science, showed that "he sensed the invisible link which connects generations, [...] the fertile strength of tradition as well as the clear foyer which has been burning for a hundred years in the house of Galois, Briot, Pasteur, and which [Tannery] has preserved with the religious fervor of a priest and the tenderness of a son. In Tannery, the École loses a living symbol of its continuity and of its ends: the visible consciousness of its nobility" [Hovelaque 1911, p. 310]. This republican "nobility" shows the model-role played by some more ancient forms of greatness [Ihl 2007]. It was often expressed in contrasting the new role played by science in ruling society with the traditional role played by religion ([Charle 1994] [Nicolet 1982, p. 310]).

Academic eulogies had traditionally played an important role in the expressions of intellectual authority. They consisted in excavating symbolically a *savant* who had died decades earlier. They were moreover the *épreuves* through which an academician showed both his eloquence and his litterary quality, thereby proving the higher merits and authority of one of the rare man of science who was able to live up to the ancient ideal of universality.[20] The illustrous ancestor was thus presented as a predecessor of the speaker, thereby linking epistemic, social and symbolic authorities through generations [Weber 2012].

---

[19] Tannery edited some of Galois's manuscripts ([Galois 1906&1907] and gave a talk at the Bourg-la-Reine celebrations [Tannery 1909].

[20] This situation can be interpreted by appealing to Max Weber's notion of specific charism, see [Weber 2012].

To be sure, Galois was nevertheless not involved in academic genealogies. On the contrary, his figure was one of continuity in the autonomy of mathematics in regard with the *Académie*. When they presented Galois as one of their precessors in the realm of "mathematical ideas," Lie and Picard thus transferred academic symbols to mathematical ones. Both Lie's and Picard's presentations of Galois may thus be analyzed by investigating how some "economies of greatness" combined elements from the worlds "inspired" and "civic" to stress a model of "merit" [Boltanski &Thévenot 1991]. This merit resorted to one's investment in the shared superior principle of pure mathematics, thereby legitimating the authority of mathematicians as a group, which indeed appeared as a sanctuary for superior minds.

Let us now consider public discourses on mathematics in broader perspective. The turn of the century was a period of growth for such discourses. In the main journal of reviews in mathematics, the *Jahrbuch über die Fortschritte der Mathematik*, the relative weight of the section "history and philosophy of mathematics" doubled from 1870 to 1900.[21] This change of quantitative scale went along with an increasing variety in the types of publications. On the one hand, some actors explicitly commented on this situation as a fragmentation of the publics of mathematics [Turner 2011]. But on the other hand, mathematicians were increasingly portraying themselves as belonging to a group at both national and international levels.

Professors played an increasing role in both mathematical institutions (such as national societies)[22] and public discourses on mathematics. The papers that were reviewed in the *Jahrbuch*'s section "history and philosophy" included reports, eulogies, obituaries, biographies, philosophical essays, etc. But the fact that this section grew faster than the whole *Jahrbuch* was partly due to the contemporary development of histories of "disciplines." These were often written by professors in connection with some boundary-works across the disciplines of teaching curriculums.

The various publications on the history of equations in the early 1890s usually presented Galois's works as an episode among others. Some attributed to Galois's achievements a specific approach to the long-term problem of the algebraic solvability of equations, i.e., no more than a landmark in a story that ended with the solutions to the general quintic given by Hermite, Kronecker, and Francesco Brioschi by appealing to elliptic functions ([Cajori 1893], [Pierpont 1895]). In this context, Jordan's *Traité* was usually presented in continuity with some previous works on substitutions [Rouse Ball 1888]. Others insisted on the "marvelous discoveries" that resulted from the impossibility to solve the quintic by radicals, even though they did not point to groups but to graphical and numerical methods for solving equations [Aubry 1894-97]. In contrast with the rupture role that would commonly be assigned to Galois later on, most presentations of the history of equations in the early 20th century were thus highlighting the continuity of some other lines of developments than Galois theory, such as modular equations, Sturm's theorem, Augustin-Louis Cauchy's substitutions etc.[23]

---

[21] On the involvement of the history of mathematics in the development of international institutions at the turn of the century (such as the French-German encyclopedia), see [Gispert 1999]

[22] See [Gispert, and Tobies 1996].

But after 1897, most historical discourses on the theory of equations would celebrate Galois's works, thereby following Lie and Picard whom they often quoted word-for-word: obituaries [Hadamard 1913a], historical accounts ([Pierpont 1897], [Boyer 1900], [Rouse Ball 1907]), the French-German encyclopedia [Ehrhardt 2007, p. 643], advanced textbooks [Weber 1895], [Burnside, 1897, p. v] as well as more elementary ones [Comberousse 1898].

Let us consider the example of the issue of who should be considered as the founder of group theory. This issue was quite popular at the turn of the century in connection with the teaching of group theory. In the 1880s, Kronecker and his follower Netto had insisted that Cauchy had played much more important a role than Galois. But after 1897, most authors followed Picard in attributing a much more far reaching concept of group to Galois "for the group concept extends far beyond [algebra] into almost all other parts of mathematics" ([Miller 1903], [Burnside 1897]).

The illocutionary power of Picard's 1897 celebration of Galois was in solidarity with his position [Bourdieu 2001, p. 159-173], or, more precisely, with his ambitioned position at this point in his career.[24] It was indeed only after World War I that, as a result of Picard's recurrent discourses,[25] Galois's works would be commonly presented as the ending of the traditional theory of equations. In the early 20$^{th}$ century, echoes of Picard's Galois were resounding in the various stratas of the public sphere, from philosophical essays to popular dailies. This situation was not unusual in regard with the way the elite of the Parisian science was expressing the symbolic prestige of its cultural authority in the medias [Broks 1996]. For the purpose of maintaining the public's conviction in the authority of the elites, the recurrent expression of the greatness of some scientific personae was at least as important as what was actually read and understood by the public, as is usually the case with the power of pageantry [Weber 2012].[26]

### 1.3. Creativity and the higher mathematical-self

The expression of a persona typically presents the qualities of an individual as exemplar of some collective values, identity issues, ethical dispositions, etc. We have seen that both the role taken on by Picard and the role assigned to Jordan regards Galois point to public involvement as a valued character of the mathematical-self at the turn of the century in France.

As a matter of fact, Picard's efforts for mediating mathematics to the public were emphasized in the daily *Le Figaro* when the latter was elected to the *Académie française* in 1924, and again at the time of his death in 1941. Several examples show that public involvement was actually one of the modern scientist's quality at a European scale, such as Mittag Leffler's political involvement [Turner 2011], Volterra's emphasis on Henri Poincaré's public activities [Volterra 1913], Émile Borel's *Revue du mois*,[27] or Gaston Darboux's insistance on the importance of keeping the door of the academy opened to the "external world" (and its generous donators) [Darboux 1911].

---

[23] See, among others, the 1907 extended French edition of Rouse Ball's history of mathematics by Montessus de Ballore.

[24] Picard's elocutionary power is also exemplified by the immediate circulation of the latters claims on Galois's achievements on elliptic functions (these claims were based on some manuscripts that would not be published before 1906) [Ehrhardt 2007, p.641].

[25] On the recurrent media depiction of Galois as a hero of French science, see [Picard, 1900, p. 63; 1902, p. 124-125; 1914, p. 98-99; 1916, p.12; 1922, p. 281-283]

[26] The notion of pageantry has been used for an analysis of the intellectual authority in [Weber 2012] by appealing to some investigations on the elaboration of the rhetoric notion of persona in the antiquity [Guerin 2007].

Moreover, public involvement was a part of the ideal of universalism of science. It was recurrently opposed to the increasing fragmentation of sciences in general, and to the increasing impenetrability of mathematics in particular. Louis de Broglie's obituary of Picard in 1942, emphasized the persona of a universal scientist in contrast with both specialized professionals and eclectic amateurs [Broglie 1942, p. 2]. In this context, the political dimension of the figure of Galois resounded differently than a few decades earlier. In the 19th century, the figure of Galois had indeed used to contrasting sharply with the image of the mathematician as devoting a plain and solitary life to science [Albrecht and Weber 2011].

Sciences were indeed especially valued during the *Belle Époque* in France. Scientists were considered as working for the well being of humankind. As is exemplified by the national celebration of Louis Pasteur in 1895 – that same year, Pasteur had also been honored altogether with Galois as a hero of the *É.N.S.* - the Third republic celebrated its *savants* by organizing national funerals [Duclert & Rasmussen 2002]. The public dimension of these celebrations thus played a key role at a time when not only had "popular science" become an object of mass consumption,[28] but when the republican ideology was also opposing science to the traditional social role of religion [Nicolet 1982] especially in organizing what Christophe Charle has designated as a cult for the masses [Charle, 1990, p. 28-35].

In this context, Lie's discourse on Galois appealed to a typical public presentation of sciences as both belonging to some outer world of ideas and in solidarity with society through their applications.[29] As a matter of fact, Lie had introduced mathematics as a model for other sciences. He had then focused on the dichotomy between Galois's "ideas" and their applications: it was the wide range of their ulterior applications to the various "branches of mathematics sciences" that proved the "unifying power" of the notion of group Galois had originally introduced in "an interesting example." It was thus the "fecundity" of Galois's ideas that proved their "fundamental" nature. Their author had thus to be honored as one of the main "creative geniuses" of the 19th century.

But while the public dimension of sciences had been usually related to their technological applications, Picard introduced some changes in the traditional hierarchies of values. He indeed attributed the prominent role to pure mathematics. In doing so, he claimed for himself the powers of definition and of representation of science that traditionally belonged to the academy. The celebration of the "heroes of the republicanized science" was indeed "both the symbol, the mean, and the end of a social stratification" [Weber 2012] in which the "authority is in the hands of the best, the most virtuous and the most meritorious" [Ihl 2007]. More precisely, official celebrations were often both characterizing species of men of science and establishing hierarchies between such species. Galois clearly characterized the "creative" specie of mathematicians in connection with the ethical dispositions of passion and engagement. Later on, the creativity of other grand mathematicians would indeed often celebrated in comparison to Galois, e.g., Poincaré [Hadamard 1913a, p. 394 & 1913b, p. 635] and Painlevé [Denjoy 1934].

---

[27] See how Borel insisted in one his reports on Cartan's application to the academy on the increasing capacities of the latter to make his ideas public [Gispert & Leloup 2009, p.51]. See also how Borel contrasted the public dimension of his book *Le hasard* with pure mathematical researches [Gispert & Leloup 2009, p.79].

[28] Cf. [Bensaude-Vincent 2003] ; [Bensaude Vincent & Rasmussen 1997].

[29] Compare to the discussion on Berthelot's reception at the *Académie française* in [Weber 2012].

Let us consider more precisely the mathematical persona on display in Picard's celebration of Galois. The latter insisted that the "general ideas" Galois had introduced in a special case of "application" (i.e., algebraic equations) highlighted the tension in "modern mathematics" between the "artistic point of view" of general theories and the "particular applications" which are most of the time necessary for a general theory to acquire a place in science [Picard 1897, p. 339]. That mathematics was not only a science but had also a cultural value in itself was a recurrent topic in public discourses. Authorities indeed claimed that mathematics, as poetry or the arts, also belonged to the spheres of feelings and imagination. In the context of the public discussions on the reform of the teaching of mathematics, Borel insisted that the issue was not to develop new points of contacts between mathematics and "*la vie moderne*" - because these contacts were already countless and would continue to increase by themselves - but to bring "*l'esprit public*" to a more precise notion of what is mathematics and what role it plays in "*la vie moderne*" [Borel 1904, p. 431-440].

In 1910, Poincaré presented a distinction between two types of mathematicians. On the one hand, some laborious carvers tackle particular cases in minute details, they aim to reach perfection through a patient analysis. On the other hand, some appeal to the divination of intuition to develop general views on vast landscapes. But these landscapes remain quite vague despite their beauty. Their authors are thus more poets than artists [Poincaré 1910, p. vii-viii].

In regard with such claims, Galois's ideas about the nature of mathematical writings could find echoes in the Romantic Movement in literature [Albrecht & Weber 2011]. In his review of the 1897 reprinting of Galois's works for the Jesuit journal *Études*, de Séguier insisted on the literary qualities of Galois who he qualified as having a "*plume d'écrivain autant que de savant.*" Moreover, Séguier claimed the esthetical value of the splendor of the truth of the "mathematical beauty" which he opposed to some utilitarian perspectives on mathematics by quoting Picard word-for-word. As for Pierpont, he emphasized how touching is the story of Galois's life by highlighting the poetic dimension of the fragments of letters that Chevalier had published in 1832.

The theme of creativity was thus tightly linked to the one of the public. For Lie, it was through applications that Jordan had mediated to society the mysterious conceptions of Galois, i.e., the "great power" of the "eternal world of ideas." But for Picard, Séguier, and Pierpont, it was more the dichotomy ideas/public that articulated mathematics with society, like with other activities culturally valued such as the arts or poetry.

The mediation between the "public" and the higher sphere of "ideas" was in a sense materialized by textbooks. Several elementary textbooks indeed presented Galois's theory as a mysterious distant horizon. Other, more advanced textbooks, such as Picard's *Traité d'Analyse* or Vogt's *Algèbre* appealed to the esthetical qualities of "Galois's theory" for legitimating a presentation of a basic version of it (i.e. limited to prime degree equations) even though it was not immediately connected to their respective topic (differential equations) or public (trainee teachers).

"Creative" works were repeatedly contrasted with "interesting" ones in reports on doctoral thesis.[30] Here, originality and creativity highlight the increasing value attributed to the notion of "research" since the 1880s. This evolution was a part of the large-scale evolution of the institutions of the teaching of sciences and of the professionalization of mathematicians as university professors For instance, in the context of the creation of new professorships at the turn of the century, Paul Appel distinguished between elementary professorships and the ones that were to favor research works in the mathematical sciences. But, Appel insisted, both were to focus on the "passionate research of truth" in contrast with the applications that were promoted in most teaching institutions in Europe and the U.S.A. [Gispert 1991, p. 63].[31]

The notion of "research" was also participating to the stratification of the emerging national mathematical communities in the "international space." In France, it especially laid the ground for a separation between the status of university professors and the one of high school teachers [Gispert 1991, p. 80]. More generally, Laura Turner has shown how "research" was a category Mittag Leffler recurrently appealed to in connection to the notions of "contribution" for establishing hierarchies not only between various stratas of practitioners of mathematics but also between nations themselves [Turner 2011]. The notion of research was thus instrumental in the construction of a specific kind of elitist and hierarchized international space.

In 1917 the Oxford professor of natural philosophy A.E.H. Love devoted his presidential address to the London mathematical society to contrast the "merit" due to the "researches of high rank in the hierarchy of works of the mind" with the "laborious" teaching issues and the "frivolous" mathematical recreations. Love appealed extensively to the figure of Galois to lay the emphasis on the esthetical nature of mathematics as a world of ideas as opposed to formulas and computations [Love 1917 p. 271]. In contrast to the articulation between mathematics and physics, novelty was presented as the primordial quality of pure mathematical research and was thus compared to creativity and the arts again.

Let us consider the additional example of Picard's address at the fifty years celebrations of the French Mathematical Society and of the French Society of Physics. "Mathematicians," Picard claimed, "do not only work for Science," they are also "artists and poets, and the word elegance is often on their lips. A geometer is not only a logician […] his reasoning needs fineness as much as order and accuracy, and without imagination, there is no spirit of invention." As shall be seen later the allusion to logicians was targeting German mathematicians. For now, let us focus on the fact that Picard insisted on the role played by imagination in mathematics by comparing Galois to Fontenelle, a figure of the Pantheon of

---

[30] Cf. [Gispert 1991, p. 79]. See among other Picard's report on Châtelet's thesis in [Gispert 1991].
[31] Cf. [Gispert and Tobies 1996].

French literature: "The seeds are often more important than the plants themselves. The art of discovery in mathematics is more precious than most of the things that are discovered."

But Picard nevertheless opposed Fontenelle's vain ambition to judge on the value of contemporary researches to the "modest sentence the great mathematician Galois wrote as a kind of will a few hours before his premature death: 'Science is the work of the human mind who is more destined to study than to know, to research than to find the truth'." What then followed was the quotation of Galois I have chosen as an epigraph for the present paper. Galois's sentence originally contrasted the *gens du monde*'s common belief in the regularity of mathematics with the quite erratic researches of the "analysts" who only find truth after having collided with walls in all directions. In Picard's discourse, this quotation supported an opposition between the mysteries of an inner circle of higher-level research-mathematicians and a public knowledge of science.[32]

### 1.4. On Galois's mathematical Frenchness

In celebrating Galois, Picard not only took his part in the shaping of a modern mathematical persona but he also stepped out of mathematics as one of the authorities of the "French science."[33] At the turn of the century in France, the celebrations of icons were indeed instrumental to the construction of a close relation between Science and the Republic. This relation participated in turn to the "moral economy" of the interactions between the masses and the elites [Weber 2012]: Science was indeed supposed to lead to wealth, supremacy, and triumph.

At a time when scientific and national geniuses had turned into switchable notions [Charle 1994], the recurrent national – and in fact nationalistic- tonalities in Picard's discourses were coherent with the latter's role as a national elite of science. After Galois had been aggrandized to the level of one who merited entry into the pantheon of Science, he became involved like other *grands savants*, in nationalistic anti-German discourse. Picard's claims on Galois's mathematical frenchness were for instance reproduced in the paper "*Les allemands et la Science*" published in 1916 by the politician Paul Deschanel, and again in 1919 in the book *La France victorieuse, paroles de guerre* which appeared a few months before its author was elected President of the French republic [Deschanel 1919, p. 122].[34]

Patriotic tonalities had been associated to the figure of Galois since the 1830s. When he first mentionned Galois's name publicly in 1843, Liouville referred to the latter as "our ingenious and unfortunate compatriot." At the time, this designation might have targeted the Italian Libri with whom Liouville was quarelling [Ehrhardt, 2010]. But this tonaly resounded again in Liouville's 1846 *avertissement*. A few years later Terquem explicitly claimed that "with Galois, France would have been granted with an Abel if a violent death had not broken the frame of such a short and turbulent life" [Terquem 1849, p. 452].

At the turn of the century, Picard repeatedly addressed the issue of the impact Galois would have had on French science if only he had lived longer. Promptly following the 1897 reprinting of Galois's works, a review of Heinrich Weber's 1895 *Lehrbuch der Algebra* highlighted how Galois had introduced the "fundamental ideas" of Algebra as it was practiced in Germany; had he lived longer, all "French Science" would have had a different orientation [d'Esclaybes 1898, p. 416]. In 1900, Picard gave a conference at Clark University (U.S.A.) on a historical presentation of *"L'idée de fonction depuis un siècle."* At this occasion, he contrasted the brilliance of the "French school" of mathematical physics with regretting the emphasis this school had laid on applied researches. According to Picard, Galois, "if he had lived longer, would have reestablished the equilibrium in pulling researches back to the highest regions of the pure theory ; the death of Galois was thus an irreparable misfortune for the French science" [Picard 1900, p. 63].

---

[32] On the function of mystery in the *Académie*'s prestige, see [Valéry 1945].

[33] Recall that it was crucial for the Parisian academic elite to demonstrate its republicanism in a context in which both institutional and political supports tended to focus on the more visible leaders of the "French science" [Fox 1984, p.120].

[34] See also [Moureu 1928, p. 334].

At the dawn of the war, while he was publicly debating the value of the "German science and philosophy" at the *Académie*, Picard appealed to the names of Cauchy, Fourier, and Galois to support the claims that the latin and the anglo-saxon civilisations produced "most of the great creations in the mathematical, physical, and natural sciences," and that "French geometers" have "opened "almost all the paths of the modern researches" in the abstract domains of pure mathematics and of mathematical physics. Galois was thus now involved in the alliances of World War I. Echoes of these alliances would resound again later in public discourses on mathematics in the 1930s [d'Ocagne 1934].

As has already been alluded to before, in 1916 Picard opposed his monograph *L'histoire des sciences et les pretentions de la science allemande* to the Manifesto of the Ninety-Three German scientists. From then on, Galois would increasingly be celebrated as an icon of the "French style of thinking" (*la pensée française*) as opposed to the German science. For instance, in 1923 the name of Galois appeared in the journal "*La pensée française, libre organe de propagation nationale et d'expansion française*." This journal was published in Strasbourg. Its frontispiece represented the French Mariane with storks wings (the symbol of Alsace), sitting on the globe of the world, with the notice "*La pensée française reigne sur le monde comme l'expression même de la liberté féconde et généreuse*." Galois was referred to at the occasion of a discussion on the masterpieces of French literature: "*Galois, mort à vingt ans, a laissé un sillage de feu dans les spéculations de la haute mathématique*" [Dunand, 1923, p. 18].

It was at the climax of the post-war propaganda on the universality of the French style of thinking that Galois's achievements came to be celebrated for their role in the "unity of mathematics." In an issue of the journal *France et monde* devoted to the topic of "the great ideas of mankind and the French style of thinking," Hadamard claimed: "Thanks to Galois's method, which might be the deepest thing a human being ever conceived in mathematics, the general problem of algebra – to which one can reduce almost everything that was studied during antiquity – is (theoretically) solved as much as it can be […]" [Hadamard 1924, p. 339].[35] In his official address for the 50th anniversary celebration of the *Société mathématique de France* in 1924, Picard laid the emphasis on the "universality of French mathematics." Alongside to Poincaré, Cauchy, Fourier, Hermite or Chasles, Galois was one of the main references Picard appealed to for supporting the inner tension of a concept which mixed the national and the universal [Picard, 1924, p. 31 -33; 1924b ; 1925].

---

[35] One year earlier, the same journal published a review on Picard's *Discours et mélanges*. This review insisted on the universality of Galois's accomplishments as opposed to the German science.

As has already been alluded to before, at the turn of the 1930s-1940s, the founders of the Bourbaki group symmetrically reversed the previous categories of the official history of Galois as an icon of mathematical Frenchness. The latter's works came to be celebrated for having paved the way to algebraic number theory as it had been developed in Germany. But the role played by the "universality" of Galois's works in the unity of mathematics nevertheless remained an idea dear to the Bourbaki group.

## 2. Broader scales of analysis and *longue durée* issues

In the previous section, we have seen how the dichotomy equations/ groups, as it had been emphasized by official discourses on Galois at the turn of the century, was related to some specific interactions between mathematic and society. We have seen also that this dichotomy pointed in turn to a more general tension between the outer world of ideas some individual heroes may access and the physical realm in which the thoughts of the "public" grazed peacefully, as André Beaunier mocked it when reviewing Galois's biography in *Le Figaro* [Beaunier 1908]. But this tension was neither specific to the *Belle Époque* nor to official discourses. In the introduction of the present paper I have emphasized that the figure of Galois escapes categories. Before analyzing further the categories used by official discourses from 1890 to 1930, it is therefore compulsory to address the issue of how specific such discourses were in regard with broader contexts and more extended time-periods.

### 2.1. Galois and Pascal: precocity and the mathematical-self in the long run

Les us first consider how Galois was referred to in the various publications that targeted the fragmented audiences of the public sphere from 1900 to 1940. We shall start our investigations with the highly representative example of a paper that was published in 1932 (i.e., at the occasion of the centenary of Galois's death) in a journal that presented itself as a "grand weekly for everybody" and which was entitled "Ric and Rac, two modern dogs." This journal was no academic publication. It stressed actualities and cultural events with pen and ink illustrations of some humorous sceneries of modern life. The paper on Galois followed a note on the political turmoil in 1832 when barricades had blossomed again in Paris and had eventually turned into a bloodbath. It was explicitly addressing the "historiographers of science" to whom it aimed at recalling a "sad fait-divers that broke the news a century ago." As was highlighted quite ironically by its concluding sentence - "by a singular whim of destiny, an eminently serious boy became the hero of a romantic drama" - the paper revolved on the tension between the variety and the triviality of the public spaces associated to the usual succession of episodes of Galois's life and the quiet place where one would expect to find "one of these severe and peaceful masters whose only duels are with equations and logarithms." After Galois had left the "world of schools," his trajectory indeed crossed barricades, a jail, the wasteland of the Glacière area, and ended at the Cochin hospital.

The tension between blackboards in classrooms and the agony in the muddy Glacière is undoubtly a dramatic one. As the novelist André Beaunier pointed it out in *Le Figaro* *littéraire* in 1908, without this tension, [36] Galois's biography would not deserve much attention, especially because Galois's story sharply contrasted with the role teachers used to devote to the exemplary lives of illustruous men. [37] Moreover, Beaunier explicitly emphasized the role played by authorities of mathematics for making Galois's story an interesting one: "Évariste Galois was a genius mathematicians. One has to know this. Further, one has to believe this for Galois's existence to be truly poignant […]. But I can nevertheless not prove Évariste Galois's mathematical genius, neither here, nor anywhere else. My reader and I shall therefore follow the testimonies the most competent people have made in memory of such an extraordinary man."[38]

The chronicler then gave a quite ironical report on the usual chain of episodes of Galois's life ("an almost virile mother," "a quite strange father," "Galois fairly understood that if he fell lonely, it was because of his genius"). But on the contrary of most public discourses that were not authored by the authorities of mathematics themselves, Beaunier nevertheless commented directly on Galois's writings:

> According to M. Adhémar, Galois's works […] came on stage only in 1870 when M. Camille Jordan published his admirable *Traité des substitutions et des equations algébriques*. It was thus forty years after Galois's tragic death that the latter's works were "discovered." These are now considered as one of the pinnacles of human thought. […] I shall nevertheless not be able to sum up Galois's mathematical ideas. I shall moreover not dare to. I shall thus quote a few sentences from the short writings left by this man of genius. These sentences are not readily intelligible for people who do not make a frequent use of mathematics. But they are nevertheless beautiful. They seem to me as a kind of prodigious hawker's pitch for selling some absoluteness.

Beaunier thus emphasized that writing on Galois involved a reflection of the writer on the legitimacy of his writings in regard with mathematics [Albrecht and Weber 2011]. After he had quoted the famous sentence in which Galois claimed that he could not carry out the process to know whether or not a given equation is soluble by radicals and that one thus had to jump over computations by considering groups of operations, Beaunier concluded:

> They use recipes! They have tricks! They are the clowns of the sublime! […] This is admirable, young, and sad! Unless I am mistaken there are more ideas here than in ten volumes of another author […] Évariste Galois, who died at twenty-one has left his ideas in the state of short and elliptic formulas. The circumstances - a little bit rough, it is true – have saved his works from a complaisant development, thereby preserving their quite strong tonality of vainglory

The fictionalization of the Galois's life had thus underlying it the identity of Galois as a genius mathematician as academic authorities had assessed it. In contrast with Beaunier who had dared commenting on Galois's writings, it was usually the attribute of precocity that embodied Galois's mathematical self. Indeed, the Ric and Rac paper was accompanied with a quite ironic pen and ink drawing of a baby, sitting next to a black board, and reading a book. As for the content of the paper, while it neither mentioned any theorem, theory, or notion, it introduced Galois as a "young mathematician," who, because he was "prodigiously gifted, was already showing when he was fifteen the promise of becoming a genius theorician, a light of the abstract science."

---

[36] In some of his unpublished drafts, Segalen had compared the mixed figure of Galois to the on of Rimbaud in connection to his essay on the *Double Rimbaud* (1906) [Ehrhardt 2007, p.3 & 665].
[37] The saillant episodes of Galois's life did not fit to traditional form of biographies of scientists [Albrecht & Weber, 2011].
[38] Compare to [Sarton 1921].

Other publications considered much more seriously the attribute of precocity, to the point of presenting it as a physical appendix. In 1900, a chronicler of the *Revue générale des sciences pures et appliquées* indeed blamed the professor Möbius, a phrenologist of Leipzig who was looking for the "bump of mathematics," for not having considered Galois's autopsy when measuring the brains of famous mathematicians such as Gauss and Dirichlet.

More importantly, the attribute of precocity pointed to the extended historical dimension the figure of Galois already had at the time. In 1833 Alfred de Vigny had indeed already compared Galois's precocity to Blaise Pascal's in a discussion on the specificities of mathematics in regard with literature [Bibas 1947, p. 272]. The association Galois-Pascal had then circulated in the long run within public discourses on mathematics. As a matter of fact, even though Liouville's 1846 edition of Galois's works has often been highlighted for its role in Galois's "resurrection" [Bourgne et Azra *in* Galois 1962] or "posthumous birth" [Ehrhardt 2010a] at the Academy, the public circulation of the figure of Galois clearly shows some continuity from the 1830s to the 1930s.

Let us now take a closer look at the early public dimension of the figure of Galois. It must first be pointed out that the question of the publicity of science was an important issue at the time. In the 1820s-1830s, the Academy was indeed forced to accept public reports on its sessions. Before the official *Comptes rendus de l'Académie des sciences de Paris* were founded in 1835, these reports were at first published in some independent journals, such as *le Globe*, *L'album national, le Temps*, or the *Annale des sciences de l'observation.*[39] Moreover, several specialized mathematical journals were founded in the first third of the 19th century. By the 1840s an emerging community of professors was claiming its autonomy in regard with traditional institutions such as the *Académie* or the *École polytechnique*.

Recall also that the world of the academic practitioners of mathematics in Paris was still a small world in the 1820s. Galois had been one of the twelve students in the science section of the *École préparatoire*. Most mathematicians at the *Académie* must have heard Galois's name and might have shared Sophie Germain's opinion on Galois as an impertinent though exceptionally gifted student [Ehrhard 2010a, p.109].[40]

René Taton has shown that Galois had been referred to publicly early on in the 1830s. He has also given a list of such public references.[41] To this list, I shall add the biographic notices on Galois that were published in Pierre-Charles Desrochers's necrology of the "outstanding people" who died in 1832 [Desrochers 1833] as well as in the successive editions of the abbot François Xavier de Feller's "universal biography" of "the men who have achieved fame through their genus, their talent, their virtues, their errors or their crimes" [Feller 1834].[42] All these public presentations laid the emphasis on the mixed figure of Galois, i.e., on the tension between Galois's passionate republican engagement and his promising mathematical capacities.[43] But even though some blamed Galois as a foolish agitator, all of them cited the latter's unpublished memoir on the theory of equations.

---

[39] This issue is highlighted in Darboux's 1901 obituary of Bertrand whose father was one of the founders of *le Globe* in 1925.

[40] A quite similar description was published in the journal *Le Précurseur* in June 1832.

[41] See [Taton 1993] for the reference to some periodicals as well as to some mentions of Galois by [Raspail 1839], Gisquet (1840), [Blanc 1842], [Nerval 1841], and Dumas (1852-1856).

[42] See also [André, 1834].

It is well known that Chevalier had published an obituary of Galois in the *Revue encyclopédique*. The latter added to his obituary some fragments of the letters Galois had sent him. Chevalier's obituary was followed by the publication of the famous testament-letter, which was originally entitled "*Travaux mathématiques d'Évariste Galois.*" Moreover, this publication was introduced with grandiloquence by the two redactors of the *Revue* who emphasized the tension between the legitimacy of the academy and the one of the public, i.e., the "knowledge of the few," as they said, versus, "everyone's rights." Even though, the authors admitted, "the knowledge of transcendental mathematics is nowadays shared only by a very few minds," they nevertheless insisted that it is the right of Galois's letter to be "contemplated by anyone who feels respect and pity." Because, they claimed, "these pages are the sacred legacy of a genius who, while he was feeling himself dying before having completed his task, turned toward humanity as if following a religious instinct, in the aim of acquitting himself by paying his tribute of new truth as a fee of passage; these are the last remains of the thoughts of a man who died before having reached his greatness."

It was thus by appealing to the human rights for pity and respect that the redactors claimed their legitimacy to celebrate Galois's "greatness." As a result, they laid the emphasis on epistemic and moral values quite different from the ones usually celebrated by academic eulogies: if Galois is unique in the "annals of science," it is, they said, because of the "strength" of the "passions" of his "ardent soul."[44] Such values were clearly the ones of the republicans and/or of Saint-Simon's followers. In his necrology, Chevalier actually insisted that even though Galois had been mostly known as an "ardent republican," he would be remembered for his scientific genius [Chevalier 1832, p. 744].

But even though the redactors claimed that through his precocious political involvement, Galois had raised the hope that the "mathematicians of the empire" might have had a successor, they nevertheless did not compare the latter to any specific great name of the mathematical pantheon. On the contrary, they presented Galois as the Thomas Chatterton of science, thereby appealing to a figure of literature especially appreciated by the Romantic Movement. It was actually in commenting this association of Galois to Chatterton that Vigny compared the precocity of Galois to the one of Pascal in 1833.

From then on, several journals in opposition to the *régime* presented Galois as a brilliant mathematician who had been a victim of the *Académie*. Much before Liouville had mentioned Galois's name again at the academy in 1843, a detailed account of the misfortunes of both Abel and Galois was given in Louis Blanc's *Revue du progrès politique social et littéraire.* The two "mathematicians" were portrayed as martyrs of the academic institution in a *brulôt* against the French educational system [Trélat 1840, p. 113]. In the same volume of this journal, one could read Blanc's call for the constitution of the *Commune*

---

[43] See [Taton, 1993] and [Ehrhardt 2007]. This lasting impact of this tension is illustrated by the well known notice on Galois Eric Temple Bell entitled "Genius and stupidity."

[44] Galois's "ardent imagination" was also emphasized in the forensic report of June 1832.

*de Paris*.⁴⁵ The fact that, much later on, the communist daily *L'Humanité* presented in 1923 Galois as a hero of the revolution of 1830 and as a precursor of the *Commune* of 1871 highlights a form of continuity in the circulation of the figure of Galois in the long run.

But in a sense, this continuous story even precedes Galois. Its main elements can already be read for instance in the editorial that opened Jacques Saigey's and François-Vincent Raspail's *Annale des sciences de l'observation.* Recall that Raspail was to meet Galois in 1830 when the latter joined the *Société des amis du peuple* (which had the former for president), and again in jail in 1831 [Raspail 1839]. Passion, youth, truth, freedom were the values of the scientific persona the *Annale*'s editors were putting to the fore. These values were associated to a presentation of the history of science that laid the emphasis on "revolutions" and quarrels of ancients versus moderns: "fertile alliances" between fields always cause resistances, Saigey and Raspail argued, "because a man prefers to break his scepter than to give it or to share it. The *savant* rarely decides by himself to unlearn. He is terrified that some innovations might humiliate and dethrone him […] and he would thus like to oppress indiscrete innovators." But the "presumptuous knowledge" of old masters is vain, the editors argued, because truth will eventually prevail with the "new generation" (which is to be sure doomed to turn into a new conservative dogma some day).

The *Annale* thus aimed at "tearing down demarcation lines" by establishing "new communication between all the parts of science" for both the sake of "the sciences" and of "the public." Progress was presented as resulting from the discoveries of links between once isolated fields of investigations, such as with the canonical example of Descartes's application of algebra to geometry and its long-term impact on analysis and mechanics. While fields of investigations are still isolated, "one is doom to walk randomly," but the "paths chance draws for us" inevitably lead to some "borders" and one thus eventually finds a neighboring science that "other intelligences had been investigating in a similar isolation." Here, one recognizes a statement very similar to the quotation of Galois chosen as an epigraph to the present paper.

Before Galois, Abel had already embodied the values Saigey and Raspail had emphasized. Indeed, after Legendre had reported on Abel's death at the *Académie*, the lengthy delay the institution had taken to recognize the value of Abel's mathematical work had crystallized public contestations. At the dawn of the 1830 revolution, Abel had been portrayed as a martyr in the journals that had forced their way in reporting on academic sessions, such as Saigey"s and Raspail's *Annale.*

Galois was already compared to Abel in 1832, e.g., in Chevalier's obituary and even before Galois's death in the paper published in *le Globe* in 1832 for supporting Galois in the trial he was facing for having given a toast for Louis Philippe's regicide.⁴⁶ Chevalier especially insisted that the similarity of Galois's and Abel's dooms was a demonstration of the "radical vice which is nowadays opposed to the progress of Science." In 1832, Saigey published a biography of Abel in which he presented the latter as an exemplary figure of "a young man who is presenting himself without any recommendation to the tribunal of science."[47] As a conclusion, Saigey urged young scholars to follow the "natural impulse" of their own "interior voice": "read and meditate the writings of men of genius; but do not become complaisant disciples or calculating admirers. Truth in the facts, liberty in the systems, shall be our motto" [Saigey 1832].

---

⁴⁵ Blanc himself mentioned Galois [Blanc, 1842].
⁴⁶The paper published in *Le globe* in 15 june 1832 is reproduced in [Taton 1982, p. 17]. See also [Dupuy, 1896, p. 198-199], [Taton 1947, p. 118].

Reflecting on the many young republicans who died when the army eventually charged the barricades in 1832, a chronicler insisted that "yes, young people are indeed precocious nowadays" for having been taught as scholars on the one hand, and for having been told not to listen to the "oppressive, tyrannical, and odious" masters of the older generation on the other hand. Fooled by the excessive value they attributed to their early successes at school, these young people were then lured to become "publicists and reformers" [Pannier 1831-1834].

The association of the attribute of precocity to Galois-Pascal exemplifies further how the long run continuity in public discourses on Galois contrasts with the discontinuity of the three main episodes of academic celebrations in 1846, 1897, and 1962. Indeed, when he edited Galois's works in 1846, Liouville clearly expressed his authority as both an academician and the editor of the *Journal de mathématiques pures et appliquées* [Ehrhardt 2010, p. 563]. In doing so, he aimed at appeasing tension in claiming he would strictly limit his comments on Galois's mathematical works. But while Liouville actually did not comment so much on Galois's mathematics, his *Avertissement* actually opposed the sterility of political passions to the fertility of experienced figures such as Richard, the professor of Galois, and the academicians Poisson and Lacroix. Moreover, Liouville incidentally disreputed Chevalier's necrology and avoided any reference to the 1832 introduction to "the testament letter" as he renamed it.[48]

But two years later, Galois's biography that was published in the *Magasin pittoresque* was nevertheless in continuity with earlier public discourses. Even though it appealed to Liouville's "high authority" to support a claim of non-partiality, the author[49] highlighted Galois's republican passion and mathematical precocity – which he compared with the one of Pascal again - by insisting on the episode of the writing of Galois's mathematical testament the night before the duel: "he seemed mostly concerned with the regret of dying without having done anything for science and for his country." The author then compared Galois's fate to the one of the "many other" young people who had died at a young age at the time, thereby implicitly referring to the bloodbath of 1832.

Encyclopedic dictionaries then usually followed the 1848 notice until the turn of the century. It is highly significant that Pierre Larousse's *Grand dictionnaire universel* attributed to Liouville the association between Galois's and Pascal's precocities as an attribute of a "superior mathematical genius," while Liouville had never made such a claim. As for the case of *Ric et Rac*, these notices always revolved around the tension between the black board and the muddy agony. The academic mathematician Joseph Bertrand blamed the inaccuracy of some of these anecdotes on Galois in 1899. But as shall be seen later, official discourses were not always perfectly accurate themselves. Claims such as Bertrand's highlight a tension between two types of legitimacies that is consubstantial to the figure of Galois.

---

[47] Caroline Ehrhardt has shown that more than a half of the memoirs submitted to the Academy were either reviewed very lately or not reviewed at all. Moreover, the first memoirs of several other young promising mathematicians such as Sturm or Liouville were rejected by the institution [Ehrhardt 2010a, p. 110].

[48] Quite similarly, the editor of Abel's collected works refuted the responsibility Saigey had laid on the academy for Abel's death (without even mentioning the name of Saigey or of his journal).

[49] The author of this biography may have been Flauguerge, one of Galois's friends, see [Dupuy 1996, p. 198.

That Galois was designated as a mathematician is exemplifying this tension. Recall that at the time, academicians were usually referred to as "geometers," and sometimes as "astronomers" in connection to the subsections of the section of mathematical sciences at the *Académie*. Other practitioners of mathematics could also be qualified as geometers by extension, even though they were usually identified by their status, i.e., by the various titles of professors, officers, civil or military engineers, priests, *élève* (pupils of the *Écoles*), *étudiants* (free students of the universities) etc. For instance, Galois was the *only* individual referred to as a mathematician in a collective volume on Paris's public institutions that was published in 1861. Here, the reference to Galois was not to be found in connection with the *Académie*, the *École polytechnique*, or the Sorbonne University, but in the entry "Cochin Hospital" where Galois had died [Audigane et al., 1861, p. 5].

To be sure, the identification of Galois as a mathematician later evolved with the complex process of professionalization of mathematicians. The emerging community of teachers of mathematics could indeed oppose the figure of Galois (as well as the one of Abel) to both the *savants* of the academy and the professors of the *École polytechnique*. That Liouville had acknowledged the quality of Richard as a teacher in his introduction of Galois's works shows that the 1846 edition had certainly aimed at appeasing oppositions between various strata of practitioners of mathematics. Later on in 1849, Terquem's necrology of Richard transferred the contestation of the academy that was traditionally attached to Galois to a charge against the *École polytechnique* whose examiners were blamed for having failed to detect Galois's genius. The episode of the rag Galois was said to have thrown at the head of one of the examiners was in continuity with the value of passion attached to the public celebration of Galois as a mathematician. Interestingly, this episode was the one Bertrand disreputed in 1899.[50]

Later on, the association between Galois's and Pascal's precocity went on circulating in public discourses even though neither Liouville, nor Picard, Lie or Dupuy mentioned it. In the early 20[th] century, it was for instance supporting a comparison between Galois and the "eleven years old apprentice pork butcher" Gioachino Rossini in a discussion on the notion of progress in music [Philippe 1913, p. 332-334]. Galois's precocity was actually explicitly discussed as an attribute of the mathematical self in connection to the various types of scientists Charles Richet aimed at characterizing in 1923 [Richet 1923, p. 104-108]. In the biography she wrote of Pierre Curie in 1924, Marie Curie insisted that Pierre Curie had rediscovered by himself Galois theory when he was a teenager [Curie 1924, p. 17]. As a character of the mathematical-self, precocity was indeed also a way to stress *a fortioti* the genius of other scientific figures, such as the botanist Noël Bernard who was also compared to both Curie and Galois [Europe nouvelle 1918].

---

[50] A quite similar discourse can be read in [Guyot 1867, p. 72]

The mathematical nature of Galois's precocity was also discussed in connection to the general notion of "creativity" and thereby in a great variety of situations which involved the specificity of mathematics in regard with other sciences or with the arts. For instance, Theodor Gomperz appealed to the precocity of Galois's "creative genius" in mathematics for validating Aristote's distinction between abstraction and induction. In reflecting on the scientific method in the case of morphology, Giard opposed the important role played by experience in sciences to the "surprising precocity of great mathematicians such as Pascal, Abel, Galois," which he connected to "the childish character" of many mathematicians and thereby to a hierarchy between sciences. This discussion was later pursued further in a comparison between "the laws of intellectual production" in literature and mathematics [Mentré 1919]. In 1913, Rougier appealed to the creative precocity of geniuses such as Galois to support the claim that teaching Latin and Greek to young students was a waste of precious time [Rougier 1913, p. 330]. In contrast, that precocity was a specificity of mathematics was also used to support the positive "influence of age on scientific personality" in other disciplines, and therefore the claim for the creation of laboratories for honorary professors [Lasseur 1933].

Appealing on the above reflection on the long-term continuity of some categories of the public discourses on Galois in the long run, one may address the reverse problem of a possible influence of such categories on some academic discourses. Official discourses such as Picard's or Lie's may indeed not only follow the codes of academic eulogies but may also be influenced by other literary forms. It is especially striking that Galois's mathematical accomplishments have been often described in the light of a "precocious" introduction of the concept of group in the history of mathematics. Hence, the role assigned to Jordan was to have made "group theory grow adult" [Dieudonné 1962].

### 2.2. The public dimension of the problem of solving equations by radicals

The finite sequence of the episodes of Galois's life is another element of long-term continuity in public discourses. The fictional nature of some of these episodes eventually prevailed despite Dupuy's attempt to establish a rigorous biography. This fictional nature was criticized by Bertrand [1899], Paul Mansion [1910], and Alain [1909] among others.[51] Mathematical issues were nevertheless spared by criticisms. Actually, the philosopher Alain even opposed Galois's contributions to "pure mathematics," which he considered as exemplar of the possibility for human beings to access a world of "pure ideas," to the moralist nature of the "edifying images" of fictional biographies.

But the dichotomy between equations and groups was no less fictional a link than others in the chain of the episodes of Galois's story. According to Pierpont, details on the episodes of Galois's life were "indispensable" to the understanding of the "early history of Galois's theory." Indeed he claimed, because "the manner in which Galois's theory of equations became public" is "intimately connected with Galois' life," it thus affords "an opportunity to give some details of his tragic destiny." [Pierpont 1899]. As said before, the novelist André Beaunier made it clear in 1908 that Galois's story would have little interest if Galois had not been a mathematician. After having presented the usual chain of biographic episodes, he thus commented on a few Galois's quotations related to the groups of equations.

---

[51] On this issue, see [Rothman 1982], [Taton 1993], [Ehrhardt 2010].

Let us consider further the example of how the dichotomy equations/groups appeared in connection to Galois in the *Revue de Paris* in 1934. It must first be pointed out that the reference to the episode equations/groups was not different in nature from the one the same journal made four years earlier to the episode of the failure at Polytechnique. In 1930, Galois had indeed been considered as a "*grande âme stendhalienne*" because he was a "*polytechnicien manqué,*" like Julien Sorel, the main character of Stendhal's *Le Rouge et le Noir* [Thibaudet 1930, p. 328]. In 1934, the reference to Galois appeared at the occasion of a review of the collective volume *L'Histoire de la Troisième République*. The chronicler marveled that the book also investigated the history of mathematics, which he described as "a thick brass wall behind which an impenetrable mystery is hidden to my sight." The academician Maurice d'Ocagne, the author of the history of mathematics section of the book, was thus in charge of opening partially the gate to the world of idea. For doing so, he stressed the importance of Galois's works by claiming (falsely) that Galois had determined all algebraic equations solvable by radicals, thus putting an end to the traditional theory of equation and opening the gates of group theory [Ocagne 1934, p. 397].

The association between Galois's works and the "theory of equations" must be considered in the perspective of the role such a theme traditionally played in public discourses on mathematics. It was indeed mostly within public discourses that Galois's works were presented in connection to a long-term history of the solvability of equations by radicals. In contrast, most mathematical papers that referred to Galois until the 1880s did not discuss the general principles of Galois"s *Mémoire* but rather pointed to other issues, usually in connection to the algebraic-arithmetic properties of elliptic functions. This attitude was characteristic of a field of research that had developed between the 1820s and the 1850s and which Catherine Goldstein and Norbert Schappacher have designated as the field of arithmetic algebraic analysis [Goldstein and Schappacher 2007a, p. 26]. Issues related to this field of research were nevertheless much more rarely mentioned in the public sphere than the solvability of equations by radicals. Galois himself had emphasized that the latter topic was a matter of gossips: "That the general equations of degree higher than four cannot be solved by radicals has become a vulgar truth […] by hearsay even though most geometers ignore the proofs given by Ruffini, Abel, etc., […] [Galois, 1962, p. 33].

Let us take a closer look at the public dimension of the issue of solvability by radicals. In the mid 1820s, the "irrefutable character" of Ruffini's proof of the impossibility to solve the general quintic by radicals was a common allegation in some biographic dictionaries.[52] A few years later, while they rarely mentioned Abel's works on elliptic integrals,[53] the public reports on academic sessions insisted on the latter's proof of the impossibility to solve the general quintic by radicals. Moreover, they echoed publicly the official judgment Legendre had made at the *Académie* in 1828 that "[Abel] can decide immediately if an equation is solvable" [L'Album national, 1829a, p. 260]. Later on, the exact same result would be repeatedly attributed to Galois as has been seen above with the example of d'Ocagne's discouse.

---

[52] See [Arnaut et al., 1825, p.300]. This notice also insists on the irrefutable character of Ruffini's proof of the impossibility to square the circle. Note that 1825 is the year when Abel published his proof on the impossibility to solve the general quintic.

[53] In his biography of Abel, Saigey had nevertheless alluded to elliptic functions as "a field outside of which neither scientific progress nor personal benefit was likely." This vague allusion nevertheless contrasted with Saigey's precise report on the issue of the solvability of the quintic.

The public dimension of Abel's premature death therefore sheds new light on the Galois's affair at the academy [Ehrhardt 2010a]. Both the confusion about what Abel had actually accomplished before his death, and the political agitation which did not end after the 1830 revolution, might have played a role in Lacroix and Poisson's decision to postpone the publication of Galois's *Mémoire*. The two academicians could especially not conclude on the originality of the *Galois criterion* in regard with Abel's works. Some of the latter's papers, fragments, and letters were still in the process of being published posthumously in *Crelle's journal.* As a matter of fact, Lacroix and Poisson's report on Galois's *Mémoire* pointed to a recently published letter of Abel to Legendre [Abel, 1830a,b], iin which the former had alluded to a statement close to the *Galois criterion*. It was on this allusion that Legendre had based his public evaluation on Abel's works. In contrast with the vague reference they made to Abel's general "exact rule" for distinguishing equations solvable by radicals, Lacroix and Poisson did not even mention the memoir in which Abel [1829] had stated a criterion that the *Galois criterion* had generalized. It seems likely that the academic report on Galois's *Mémoire* was partly reacting to a public interpretation of Abel's achievements on equations.

The field of reserch of arithmetic-algebraic-analysis was nevertheless in the process of emerging in the 1820s, even though mostly outside of France. The context of the emergence of this field is illustrated by some aspects of both Abel's and Galois's investigations on the general relations between the roots of algebraic equations. These investigations were indeed modeled on the special equations related to cyclic, elliptic and hyperelliptic functions [Goldstein and Schappacher, 2007a, p. 34]. For instance Abel stated that prime degree irreducible equations whose roots are rational functions of one of them are solvable by radicals. This criterion was a generalisation of the properties of the roots of the division equation of the lemniscate (i.e. an elliptic function), which investigations had been modeled in turn on Gauss's treatment of the binomial equation of the division of the circle [Abel, 1829, p. 660][54]. In 1830, Galois not only generalised Abel's criterion to equations whose roots are rational functions of two of them, but he also investigated the modular equations of elliptic functions as well as the primitive equations of prime power degree on the model of Gauss's treatment of the binomial equation.

But this articulation between special model cases and general properties did not fit to the traditional definition of the theory of equations. This definition indeed revolved on a tension between the general and the numerical.[55] Moreover, as has been shown by Caroline

---

[54] The interconnection of Galois's and Abel's works on equations and elliptic functions was still emphasized at the turn of the century, see [Pierpont, 1899, p. 299].

[55] In contrast with Galois's approach to the general/the special, the topics on which Galois published (general resolution, numerical equations, continuous fractions) were fitting nicely to the definition Lagrange had given of algebra as a science partitioned in three main sections. First, the general theory of equations concerned properties common to all equations (such as the number of roots in function of the degree). Second, the *general* resolution of equations consisted in expressing

Ehrhardt, the few applications Galois presented of his general principles as well as the incompatibility of these principles with effective computations was hardly fitting the usual dichotomy theory/ applications in connection with the mathematical training at the *École polytechnique* [Ehrhardt 2010a, p.104]. In reaction to Galois's claim that his *Mémoire* was an "application" of a "general theory," Poisson and Lacroix invited the author to publish his theory as a whole. Moreover, they criticized Galois's *criterion* for providing conditions between the roots themselves instead of appealing to conditions "exterior" to the equation considered (i.e. on the coefficients) [Lacroix and Poisson, 1831, p. 660].

The evolutions of the theory of equations in the teaching of mathematics in France from the 1820s to the 1840s have been analyzed in [Ehrhardt, 2007, p. 211-236]. Liouville's edition of Galois's works was indeed contemporary to the emergence of the *Algèbre supérieure* as an intermediate discipline between elementary arithmetic and algebra and the "higher" domain of analysis as it was taught at the *École polytechnique*. Recall that despite the fact that Galois's works were scarcely taken into account in the first edition of Serret's book in 1849, the figure of Galois was nevertheless already celebrated in the framework of a long-term history of the "theory of equations." While the crucial role assigned to the issue of solvability by radicals constrasted with contemporary researches on equations, this historical perspective clearly aimed at delimiting the domain of algebra as it was taught.

For most of the 19th century, the theory of equations was actually no more an autonomous domain of research in France than algebra itself was an object-oriented discipline shared by a community of specialists [Brechenmacher and Ehrhardt 2010]. The architecture and the content of Serret's *Algèbre* were indeed implicitly built on the higher point of view of analysis. The issues related to Galois (number-theoretic imaginaries, numbers of values of functions, substitutions, solvable equations) were first steps toward investigations of the algebraic-arithmetic properties of elliptic (or abelian) functions. This "higher point of view" involved complex analysis, arithmetic considerations on congruences or quadratic forms, the algebraic theory of invariants, etc.

Let us consider the example of Charles Hermite's two first papers. While he was still in high school, Hermite had published in 1842 in the *Nouvelles annales de mathématiques* a paper on the impossibility of solving the quintic by radicals. This paper was followed in 1844 by a memoir on the division of the periods of abelian functions, which was announced at the *Académie* and was published in *Liouville's journal*. The paper on the quintic was introduced with enthusiasm by Terquem who related it to both the classification of transcendental functions into species pursued by Liouville in 1837 and to Wantzel's 1837 proof of impossibility of the antique problems of the duplication of the cube and of the trisection of an angle.[56] In doing so, Terquem explicitly alluded to the higher point of view on equations provided by the investigations on the nature of what he designated as the "*irrationnelles.*"

---

the roots as functions of the coefficients, a topic that was not limited to the solvability by radicals as can be illustrated by the resolution of the general cubic by circular functions. Finally, the third section concerned the approximation of the roots of numerical equations.

[56] As is exemplified by Liouville's works, the development of the theory of abelian functions in the 1830s had implied more and more concerns for the distinction between the algebraic and the transcendental as regard to methods, functions, and quantities.

This designation pointed to both the quantities that were defined by algebraic equations and to those that were associated to the transcendental functions investigated by Liouville in Abel's legacy. Moreover, Terquem called for the development of studies on such *irrationnelles*. The two main authors who responded to this call, Wantzel and Victor-Amédée Lebesgue, were also the first to greet Liouville's project to edit Galois's works.

Hermite's second paper was strongly supported by Liouville, to the point that it triggered a controversy with Libri who had self proclaimed himself as the specialist of Abel's works.[57] It was at this occasion that Liouville announced publicly his aim to publish Galois's works [Ehrhardt 2010a].

In his 1846 *avertissement*, Liouville nevertheless left a veil of mystery on the content of the "general method" for which Galois deserved to be considered as one of the rare true "inventors" in mathematics. Moreover, he did not say a word on the connections between equations and elliptic functions and actually did not publish Galois's fragments on abelian functions.

One may wonder about Liouville's motivations. The activity of collecting works and editing them was publicly valued at the time. Biographic notices on Liouville in dictionaries were indeed systematically listing precisely the various "esteemed editions" the latter had given.[58] In contrast, the allusions to Liouville's mathematical works were usually quite remote.[59] Moreover, it was in connection to the edition of Galois's works that Terquem explicitly spared Liouville from his attacks against the *École polytechnique* in 1849. This edition might thus have partly aimed at establishing the latter's authority on the emerging community of teachers. Although there is no evidence that Liouville's filtered discourse on Galois was intended at targeting a broad audience of practitioners of mathematics, it has been seen above that the main actual impact of the 1846 *Avertissement* on public discourses would be to substitute for decades the exact statement of the *Galois criterion* to the earlier vague allusions on Galois's marvelous achievements on equations.

Later on, while Serret's 1866 textbook commented on Galois in the framework of the theory of equations, Jordan's 1870 treaty inscribed the *Mémoire* in a general "*théorie des irrationnelles*." Unlike the public dimensions of the issue of the solvability of equations by radicals, on the more local levels of mathematical papers, Galois's works had indeed not been much commented in connection to "algebra" or to the general theory of equations. Betti, Kronecker, Hermite, and Klein all insisted that the traditional problem of expressing the roots of equations by algebraic functions of the coefficients had to be replaced by the characterization of some "orders of irrationalities." The types of quantities that correspond to types of non solvable general equations had to be characterized by some adequate analytic expressions, as in Hermite's solution to the general quintic through elliptic functions [Goldstein 2011].

---

[57] Before Cauchy returned from exile, Libri and Liouville had been among the rare Parisian academician to deal with the works of Abel and Jacobi. This involved the consideration of "roots foreign to the equations" and even of "groups of roots", see [Liouville, 1843].

[58] In addition to Galois's works, Liouville published Monge's *Géométrie*, and Navier's *Leçons.*

[59] When there was no obituary to refer to (as was the case for living beings), notices on scholars were usually based on compilations of bibliographical references from which they copied lists of publication (such as *La littérature française contemporaine*). But for specialized papers, the statement was usually reduced to something like "M. Liouville is the author of a number of important discoveries presented in a sequence of *Notes* and *Mémoires* whose titles cannot be listed here."

In Bertrand's 1867 report on the progress of mathematics,[60] Galois' works were clearly inscribed in a twofold collective dimension [Bertrand 1867, p. 3-17]. First, the "highest progress in Analysis since the time of Lagrange and Laplace" were explicitely presented as a matter concerning the "scientific superiority of France." Here, Bertand celebrated Hermite's works on the nature of the algebraic and transcendental irrationals defined by elliptic and abelian functions and their related special equations. But he also highlighted Jordan's recent works on the "groups Galois has attached to equations."

Second, the *Algèbre supérieure* "plays a key role in the researches presented above," because, Bertrand insisted, "M. Serret's book is not a textbook of algebra; the reader must be already familiar with the general and classical methods. M. Serret then leads [the reader] by an easy path to the highest results of this branch of science [...] with both precision and clarity, [the book] presents in a homogeneous unity the beautiful works of Galois, and the ones of M. Hermite, Kronecker and Betti" [Bertrand 1867, p. 14-15]. A later echo of the two fold collective dimension of Galois's works would resound in 1898 when Paul Tannery's *History of sciences in Europe* would discuss the works of Galois in connection to the approaches of Grassmann and Hamilton on generalized numbers [Tannery 1898, p. 739]. Even later in 1942, Broglie's obituary of Picard mentioned that when he was a student, the latter was "much interested in Algebra and had thus certainly already sensed the early symptoms of his vocation as an Analyst." [Broglie 1942, p. 6].

### 3. Nations and disciplines

In the 1860s, there was thus already a twofold collective interpretation of the works of Galois. For this reason, Jordan's claim to deliver the commentaries on Galois that Liouville had promised is a telling illustration of the different autonomous developments that would tear apart the field of arithmetic algebraic analysis [Goldstein and Schappacher 2007b, p. 97]. In the name of Galois, Jordan did indeed reorganize various results. Previous works that had made precise references to Galois, such as those of Hermite, thus fell into the global legacy of Galois, in the company of works that used to be disconnected from any reference to Galois, such as Cauchy's substitutions or Clebsch's geometrical problems of contacts.

A few years later, Klein interpreted Galois's works in the context of geometric transformations and invariants. He did not appeal to Jordan's approach but followed Hermite and Kronecker in focusing on the special groups attached to the modular equations of elliptic functions of order 5, 7, and 11, i.e., "the three Galois groups." Because these groups appeared at the core of Klein's various interpretations of the nature of the irrationality of the general quintic, the latter celebrated Galois for the introduction of a fully general notion of group in the special case of the theory of equations.

In the 1880s, the question of the status of the notion of group with regard to arithmetic, algebra, and analysis was much debated. Kronecker especially developed a constructive arithmetic theory of irrational quantities. Even though the latter had rejected Galois's approach, his theory would play a key role in most presentations of "Galois theory of general equations" until Hilbert and Weber would lay the emphasis on Dedekind's conceptual approach in the mid 1890s.

---

[60] This report was a part of the collection of reports to the government on the "progress of letters and sciences in France."

The disciplinary issues related to the nature of Galois groups would often have national, and in fact nationalistic, overtones at the turn of the century. Much attention has already been devoted to the inscription of Galois's theory in the emerging autonomous disciplinary framework of "Arithmetic and Algebra" in Germany.[61] In the present paper, I shall thus focus on the public claims on the power of unification of analysis in France.

### 3.1. France, analysis, and the unity of mathematics

The dichotomy between elementary equations and the higher point of view of groups at the turn of the century can be considered in continuity with the two-fold interpretation of Galois's works as presented in the previous section of this paper. We have seen also that this dichotomy reflected the duality between the public and the world of ideas, as mediated by Jordan. But in addition to his clarification of Galois theory, the role of the researcher who closed the algebraic issue of the solvability of equations was also assigned to Jordan. Lie and Picard indeed both claimed that, unlike the previous works of Lagrange and Cauchy, Galois groups had exceeded the boundaries of algebra in introducing ideas lying at the roots of the branches of modern science and whose "far reaching impact appears to us more and more every day" [Fehr, 1897, p.756]. The seeds Galois had sown in the special case of equations were to blossom into a general notion of analysis.

This claim should nevertheless not only be considered as having solely aimed at promoting Picard's or Lie's contributions to continuous group theory and differential equations. We have seen in the first section that Picard, in particular, clearly took on the role of an official public authority on mathematics. Several authorities such as Tannery, Picard, Appell, or Poincaré indeed contrasted the "richness" of the power of unification of analysis with the "poverty" of considering algebra and/or arithmetic as autonomous disciplines.[62] These official lines of discourse usually pointed to recent developments in Germany in the legacies of Kronecker or Dedekind.

Here analysis took a broader meaning than the one it would take later on as a specific discipline centered on the notion of function. In the first third of the 20$^{th}$ century, mathematicians such as Hadamard, Borel, Lebesgue, and Baire, both promoted the evolution of analysis toward an object-centered discipline (especially through the series of monograph of "Borel's collection") and perpetuated some earlier discourses on the unification of mathematics through analysis.

Before the blossoming of the self-designated "brilliant French school of real function theory," Picard was already claiming in 1890 that the notion of function is the central notion of mathematics [Picard 1890a]. Recall that, in France, the mathematical sciences were mainly divided at the time between analysis, geometry and applications. In this context, and in contrast with some contemporary developments in Germany, algebra and arithmetic were traditionally interconnected in France at an elementary level while, from a higher perspective, algebra was considered to belong to analysis. Recall that Jordan's *Traité* itself had presented Galois theory as a general theory aiming at providing a higher point of view on the classifications and transformations of the "*irrationnelles*." The general part of this theory was explicitly considered as belonging to the analysis, it was then to be applied to algebra, geometry and transcendental functions.

---

[61] See [Kiernan 1871], [Corry 1996].

[62] Some authorities, such as Koenig, contested the increasing importance attributed to analysis in France at the turn of the century. But these oppositions laid the emphasis on the traditional importance of geometry without challenging the fact that algebra and arithmetic belonged to analysis. See [Gispert 1991, p. 94]

Later on, while Drach's 1895 *Algèbre supérieure* had appealed to Kronecker for presenting the "famous theory created by Galois" as an extension of arithmetic, Tannery's *Préface* explicitly recalled that it was analysis that provided a higher point of view on "the general irrationality" of which "the algebraic number is nothing more than a particular case" [Borel et Drach 1895 p. iv]. In his *Traité d'analyse*, Picard made it clear that the algebraic Galois theory was a first step toward the higher point of view of analysis. Even later, in 1913, George Humbert's lectures on the theory of substitutions at the *Collège de France* insisted along the lines of Jordan's *Traité* on connections between group theory and elliptic or abelian functions.

As is exemplified by Picard's 1890 obituary of George Halphen, being an "algebraist" was not considered as a specialty but as a specific attitude toward mathematics. Picard indeed introduced the obituary by claiming that one can nowadays distinguish between two different orientations in the mathematical thought ("*la pensée mathématique*"):

> The ones aim before all at extending the domain of knowledge. Without always caring much about the difficulties they leave behind them, they do not fear to go forward and always look for new fields of investigations. The others prefer to stay in a domain of already developed notions which they seek to deepen further; they want to exhaust all consequences and they try to highlight the true grounds of the solution of each question. These two directions in the mathematical thought can be seen in all the branches of this Science [...] the first one can nevertheless be found more often in connection with integral calculus and functions theory, and the second one in connection to modern algebra and analytic geometry. Halphen's works were mostly related to the second orientation; this profound mathematician was before all an algebraist.[63] [Picard 1890b, p.489-490]

To be qualified by Picard as an algebraist was thus a mixed blessing. That it was less valued than being an analyst is obvious. Analysis was indeed on the side of the increasingly valued "creativity." Moreover, we have seen how Picard's dichotomy would be echoed later on by Poincaré's distinction between artists and poets. Even later, one can read in Picard's obituary of Jordan in 1922 an implicit criticism of the latter's tendency to develop very general approach to mathematical questions "as if he feared that some particularity may impeach him to see the true reasons of things. Jordan has really been a great algebraist; the fundamental notions he introduced in analysis will save his name from oblivion" [Picard 1922a, p. 210].

---

[63] Il semble que l'on puisse aujourd'hui distinguer, chez les mathématiciens, deux tendances d'esprit différentes. Les uns se préoccupent principalement d'élargir le champ des notions connues ; sans se soucier toujours des difficultés qu'ils laissent derrière eux, ils ne craignent pas d'aller en avant et recherchent de nouveaux sujets d'études. Les autres préfèrent rester, pour l'approfondir d'avantage, dans le domaine de notions mieux élaborées ; ils veulent en épuiser les conséquences, et s'efforcent de mettre en évidence dans la solution de chaque question les véritables éléments dont elle dépend. Ces deux directions de la pensée mathématique s'observent dans les différentes branches de la Science : on peut dire toutefois d'une manière générale, que la première tendance se rencontre le plus souvent dans les travaux qui touchent au Calcul intégral et à la théorie des fonctions : les travaux d'Algèbre moderne et de Géométrie analytique relèvent surtout de la seconde. C'est à celle-ci que se rattache principalement l'œuvre d'Halphen : ce profond mathématicien fut avant tout un algébriste.

Halphen's obituary also exemplifies that the richness of Analysis was clearly opposed to the works on the foundations of mathematics and thereby to both the phenomenon of arithmetization of analysis and to the consideration of arithmetic and algebra as autonomous disciplines.[64] In this context, Picard's recurrent references to "Galois's ideas" in relations to the notions of "groups" pointed to the unifying role he assigned to Analysis. The theme of the unity that analysis provided to mathematics was later mixed with the one of the universality of the French style of thinking. The idea that Galois had terminated classical algebra in introducing the general concept of group was instrumental to the claim of the essential role France had played for the alleged achievement of the unity of mathematical sciences. During the inter-war period, this claim was recurrently emphasized in the public discourses of authorities. Recall how Hadamard [1923] and d'Ocagne [1930-1936] followed Picard in attributing to Galois's ideas the role of a revolution in the whole history of mankind.

Opening the entry "Mathematics" in the 1928 dictionary *Larousse mensuel illustré* one reads the statement that "we are especially pleased to highlight the predominant role played by France in the evolution [leading to the unity of the mathematical sciences]; such is by no way a French opinion, it is an objective fact and recognized as such at an international level." This claim was accompanied by a portrait gallery presenting a few foreign mathematicians faithful to the French influence (i.e., two Italians, one Swedish and one Polish) surrounded by seven French mathematicians lead by Picard.

### 3.2. Epistemological discourses as mediations

Official discourses were more or less faithfully reported in the medias.[65] Entries in dictionaries were often either written by scientists themselves or extracted from the publication of an academic authority (such as d'Ocagne in the case of Galois in the 1930s). In the *Revue générale des sciences pures et appliquées,* official discourses - such as Picard's views on functions - were often reproduced word-for-word. But on the other hand the latter journal was also offering a forum for various practitioners of sciences including the ones who, like the polytechnician engineer Autonne and the jesuit abbott Séguier, were enthusiastic about the emphasis some German mathematicans laid on algebra and number theory. One could thus read in the book reviews section a praise of Séguier as one of the "most eminent contemporary algebraist" [Autonne 1913].

In some other journals, chroniclers were mediating official discourses. Some would abstract them, some would reduce them to a selection of citation, some would settle for some mere allusion, and some would reformulate them in developing new original presentations. In any case, the editorial orientations played the key role. Picard's views were for instance always celebrated in the catholic *Revue pratique d'apologétique*, they were reported on with much respect in *Le Figaro* while in contrast *L'humanité* usually mocked or criticized official discourses of the authorities of science after 1916.

---

[64] On the arithmetization of mathematics, see [Petri and Schappacher, 2007], and [Corry, 1996].

[65] In the case of mathematics, it does not seem that the publicists of the journals of "popular sciences" played a similar ambivalent role as the one they played in regard with discourses of officials of other sciences [Bensaude Vincent 2003, p. 151-154; Saint-Martin 2008]. The notion of "popular mathematics" should nevertheless be studied further.

Some papers were in between reports on official discourses and original contributions. These appealed to other types of legitimacies than mathematics, i.e., philosophy, literature, or politics. We have seen how a critic and novelist such as Beaunier dared to comment with irony on Galois's mathematical writings. Let us now investigate a coherent corpus of papers that aimed at presenting some recent mathematical works through a "philosophical march" [Winter 1908, p. 323], thereby acculturating these works to the literary form of the shared culture of the elite. These were mostly published in some periodicals such as the *Revue de métaphysique et de morale*, the *Bulletin de la société française de philosopie*, but also in the series of monograph of the *Nouvelle Collection Scientifique* directed by Borel at the *Librairie Felix Alcan.*

To be sure, this corpus deserves to be investigated for its own sake.[66] But here we shall focus on questionning these epistemological essays for their mediating role in public discourses on mathematices. Scientists and mathematicians were among the authors of such essays. Through Borel's *La revue du mois,* their views were echoed in the general press. Beaunier's papers on Abel and Galois were for instance both triggered by two papers published in *La revue du mois* by the the mathematicians Mittag-Leffler and Adhémar respectively. But we shall see that some authors also appeal to the legitimacy of philosophy for developing their own views, thereby resisting in a way to the public authority of mathematics.

The corpus of French philosophical papers that debated the question of the relationship between arithmetic, analysis, and algebra has a strong intertextual coherence. One of the main shared reference is a paper published by Louis Couturat in 1898 in the context of the controversy between Poincaré and Bertrand Russel on the notions of numbers, magnitudes, and on the foundations of geometry.[67] At this occasion, Couturat opposed the "science of order" to the process of arithmetisation of mathematics he associated to Kronecker. Couturat especially blamed the autonomy attributed to "pure algebra" and argued against the central role devoted to the notion of integer in Kronecker's arithmetic theory of algebraic magnitudes. The aim of grounding mathematics on the notion of integer, he said, amounts to reducing mathematics to its "poorest part,"[68] because the "idea of order" is irreducible to the notion of number. Couturat then supported his claims by the example of the transversal role played by group theory in analysis, arithmetic and algebra.[69]

The science of order was presented as a French tradition which originated with Descartes and was developed later by Poinsot and Galois. The latter's merit, Couturat claimed, "was precisely to have perceived the use algebra could make from the notion of order even though the latter is apparently foreign to algebraic speculations" [Couturat 1898, p. 436-440 & 445-447]. This claim echoed Picard's and Lie's presentations of Galois's works. The science of order was described as transversal to algebra, arithmetic and mechanics in the sense that it focused more on relations than on objects. It especially incorporated group theory and topology. Before Couturat, Jordan himself had presented his works as belonging to the "theory of order" in his thesis in 1860 and again when applying to the Academy in 1881 [Brechenmacher 2012].

---

[66] On the French epistemological milieu at the turn of the century, see [Castelli Gattinara 2001], [Vogt 1982], [Nye 1979].
[67] Cf. [Nabonnand 2000].
[68] The definition of the concept of number was characterised as an "illogisme déconcertant" in [Dufumier, 1911, p. 729].
[69] On Poinsot's theory of order, see [Boucard, 2011].

Couturat's claims were somehow contradicted in two papers published in 1908-1910 by Maximilien Winter, one of the actors involved in the foundation of the *Revue de métaphysique et de morale*. The first paper, "Sur l'importance philosophique de la théorie des nombres," aimed at characterising as well as evaluating the relations and relative roles of the "three levels of pure mathematics," i.e., "analysis, algebra, and arithmetic" [Winter 1908, p. 321]. Even though Winter acknowledged the intimate relationships between the three levels, he nevertheless praised the autonomous value of both "number theory" [Winter, 1908] and "modern algebra" [Winter, 1910].

Winter indeed claimed that philosophers had some legitimacy for commenting on mathematics. More precisely, he contested the idea attached to Auguste Comte's positivism that "philosophers" should only discuss the "general results" of sciences [Winter 1908, p. 323]. On the opposite, he argued, philosophers have to deal with the "technical forms" of sciences for developping some external evaluations of sciences: "such a critique is to sciences what drama critique is to theatre plays." Such a claim gives a fine characterization of a type of papers whose purpose was to mediate some recent scientific works to an elite public which in turn gave these publications their legitimacy.

For instance, Winter claimed that philosophers had to deal freely with historical orders of developments in the aim of highlighting the genesis and filiations of the "deepest notions," as opposed to the "fastidious chronologies" which make the history of mathematics look like a "phonebook." Such a tension between two types of legitimacies also involved issues of litterary style. In contrast with the very simple style - sometimes reduced to an accumulation of facts – of most of the historical perspectives that were developed in connection with the teaching of mathematics, papers such as Winter's presented a litterary quality similar to the one of official academic discourses.

Winter's role as a mediator is illustrated by the mix nature of the references the latter appealed to. These involved official discourses of authorities such as Picard, academic historical works such as Paul Tannery's, both French and German treatises or textbooks such as [Cahen 1900] and [Weber 1895], and Georges Humbert's contemporary lectures at the *Collège de France*. The public dimension of this mediation was crucial. Indeed, Winter recognized that his discourse had to remain on an "elementary level" to reach its audience. Moreover, his philosophical critique explicitly aimed at contributing to the "elementarization" he considered as necessary for the progress of science. Winter actually militated for the development of an "algebraic and arithmetic culture" by inscribing these two domains in the programs of the Faculties of sciences. Such an ideal was hardly compatible with the tendency of some of other types of public discourses to lay veils of mysteries on the marvels of Galois's mathematical achievements, such as when Hadamard insisted on the "marvelous capacity" of Galois's method to bridge ancient mathematics with new mathematics or when Couturat presented Lie's theory of transformation groups as "a science and a method of marvelous fecundity."

In 1910, the essay "*Caractères de l'algèbre moderne*" was grounded on Weber's presentation of Galois theory, and thereby on Dedekind's notion of *Körper*. From this retrospective point of view, Winter developed a nuanced presentation of the complex history of the various forms of references to Galois in connection to the problem of the "irrationals" in the works of Hermite, Kronecker, Brioschi, Jordan, Klein, and Dedekind. The essay thus presented an original synthesis between the non-disciplinary transversal perspectives of elliptic functions and the object-oriented perspective of algebraic number theory. It was actually on this epistemological synthesis that Winter grounded his claim for the legitimacy to consider algebra an arithmetic at the same level as analysis [Winter 1910, p. 497&528].

The opposition between the richness of Analysis and the poverty of Algebra was nevertheless emphasized once again by Hadamard in 1912 in a paper that stressed the central role played by the notion of function in the long run development of mathematics. That this notion was already "hidden" in antique mathematics (such as astronomical tables, magnitudes etc.) was one of the arguments put to the fore to express the essential nature of functions as opposed to Kronecker's focus on the concept of number [Hadamard 1912].

Later on in 1920, Pierre Boutroux's *L'idéal scientifique des mathématiciens dans l'antiquité et les temps modernes*, explicitly appealed to the legitimacy of mathematicians as opposed to the one of philosophers: "this is a book written by mathematicians for mathematicians" [Boutroux 1920]. Following [Picard 1900] and [Hadamard 1912], "Modern algebra" here designated Descartes's "method of combinations." It thus pointed to a historical episode that was to be followed by the unfolding of the concept of function and the development of analysis. As Léon Brunschwicq summed it up, "M. Boutroux's central idea is the following : with Descartes's algebra and Leibniz's differential calculus, mathematics has taken on a specific form which has become classical and which has even been extended recently with the arithmetization of analysis and with the logistician panlogicism; but this specific form is nevertheless no more than an episode in the evolution of mathematics, an episode whose greatness and decadence may both be analyzed" [Brunschvicg 1923].[70] Galois theory was presented as one of the main proof of the failure of algebra and logic [Boutroux 1913 p. 129].[71]

Boutroux presented mathematics as an eternal duel between the human mind and a rebel matter whose nature is essentially transcendent and thus escapes both logicism, nominalism, and conventionalism. He thus especially highlighted Galois's presentation of the progress of analysis as a random walk, which he quoted in full length. The citation in the epigraph of the present paper was thus opposed to the superficiality of considering algebra and logic as self-sufficient languages in which one can build objects of investigations by systematic combinations. "There is no Berlitz school for mathematical analysis," claimed a reviewer of Boutroux's book in the *Revue de métaphysique et de morale* [1921, p. 6-7].

---

[70] L'idée central de M. Pierre Boutroux est celle-ci : la forme que la mathématique a revêtue avec l'algèbre de Descartes, avec le calcul différentiel de Leibniz, qui a tendu à devenir classique, qui s'est même prolongée de nos jours par l'arithmétisation de l'analyse et par le panlogisme des logisticiens, ne représente qu'une période dans l'évolution des mathématiques, période dont il est permis de suivre la grandeur et la décadence.
[71] In this context, Galois was thus associated to Descartes and Leibniz. See [Boutroux, 1919b].

As for Winter, even though his review on Boutroux's book contested the latter's simplified presentation of algebra [Winter 1919, p. 665], he nevertheless recognized the essential nature of the notion of function by acknowledging Hadamard's authority whose 1912 paper he quoted word-for-word.

### 3.3. Germany, algebra, and the unity of mathematics

In 1923, a paper of the communist daily *L'Humanité* reported on a group of pupils of the *École Normale Supérieure* who had commemorated the 1871 *Commune de Paris*. The students were presented as followers of the revolutionary Galois whom the paper opposed to the "the French official science" [L'Humanité 1923, p.4]. The paper on Galois was published side by side with a pen and ink drawing of Jean Jaurès - the mythical socialist leader and founder of *L'Humanité* - that illustrated a paper on the latter's "intellectual life."

The 1923 paper also celebrated the "fighter Galois," thereby implicitly referring to the memory of the war. It was indeed after the peak of the war (and the Russian revolution) that *L'Humanité* increasingly criticized the official French science.[72] In the 1920s Painlevé would for instance be recurrently blamed as a war criminal for his governmental involvement as a mathematician. But recall that 1923 was also the year when a monument to the dead was erected at the *É.N.S.,* as one of the various expressions of the cult of the dead that followed World War I in France [Beaulieu 2009, p. 13]. The flippant attitude several *normaliens* adopted in regard with the teaching at the *É.N.S* and, more generally, with traditions [Goldstein 2009, p. 166] certainly gave a new actuality to the rebellious facet of the figure of Galois.

In 1922-1923, the Bourbaki-to-be André Weil, Jean Delsarte, Henri Cartan, Jean Dieudonné and René de Possel were beginning their studies in mathematics. When recalling their training years, these mathematicians pointed to the space taken by the ghosts of the *É.N.S.*, i.e., those who had died at the war before they had a chance to prove their mathematical quality, thereby abandoning their successors by their own [Beaulieu 2009; Mazliak forthcoming]. But even though the Bourbaki would claim they were deprived of any model, of any "master",[73] it seems as if the figure of Galois had provided the model of a ghost. Recall that already in the 1830s, the latter had been described as the mathematician who had died before he had a chance to prove his grandeur. The topic of what Galois would have given to both mathematics and to France if he had lived longer had been discussed recurrently for a century. In the 1920s, this issue echoed the hypothetical history of what the world had missed with the missings of the war.

As has been alluded to in the beginning of this paper, the post-war generation of *normaliens* who visited Noether and Artin in Göttingen and Hambourg in the late 1920s symmetrically reversed the previous categories of the official history of Galois as an icon of mathematical Frenchness. The latter's works came to be celebrated for having paved the way to algebraic number theory as it had been developed in Germany. Should this evolution be partly considered as a form of reaction against the ideology of the universality of the French style of thinking as heralded by some authorities of French mathematics? It is well known that the Bourbaki strongly blamed several authorities from the older generations, i.e., the "pundits of the French science" as Weil called them [Gispert and Leloup 2009, p.41]. Moreover, they refused to follow the earlier models of mathematical life which emphasized patriotism, discipline and the patriarchal family [Goldstein 2009, p. 174]. From the cynical individual attitude of the ones to the political engagement of the others, the figure of Galois could serve as a banner for a variety of conceptions on collective involvements: we have seen that this figure had indeed been traditionally opposed to the model of the "*grands hommes*" the younger generations were urged to follow after the war.[74]

---
[72] In 1916, *L'Humanité* celebrated the publication of the collective volume *La science française*, see [Lafitte 1916].
[73] This expression comes from. André Weil. Cf. [Goldstein 2009].

The fact that the Bourbaki did not even allude to Picard's public discourses on Galois's achievements highlights how former models of actions were rejected. This situation is exemplified by the public dimension that would be given to the opposition between the archaic presentation of Galois's theory in Serret's *Algèbre* and the "modernity" of van der Waerden's *Moderne algebra*. But even though Serret's textbook was indeed still edited in the 1920s, several updated presentations of Galois theory had nevertheless been published in France since the 1890s. Moreover, some prominent German textbooks – such as Weber's – had been translated in French. The tension between discourses and reality here is very similar to the one that has been investigated by Catherine Goldstein for the case of number theory. The war had indeed caused important ruptures in both memories and intellectual traditions [Goldstein 2009]. But the situation that resulted was nevertheless complex and certainly not limited to the importation of new field of research from Germany because of a generational gap.

As has been shown by Juliette Leloup, the doctoral thesis that were defended during the interwar period in France depict a mathematical landscape quite different from the one that has been emphasized by retrospective memorial discourses. In contrast with most actors' testimonies, France was indeed not a desert in algebra and number theory in the early 20th century and some actors of the older generations were actually bridging the pre-war and post-war eras [Goldstein 2009]. Moreover, the first Bourbaki were still anchored in the tradition of French analysis [Beaulieu 1993], and the works of some of the key actors of the pre-war generation such as Borel or Lebesgues could not be reduced to their past involvement in the "brilliant French school of function theory" of the turn of the century [Gispert & Leloup 2009].

The investigations of the present paper on some public expressions of mathematics highlight how the categories used by some discourses on the interwar period have underlying them some long-term continuities despite the ruptures these discourses emphasized. Let us consider the example of "modern algebra." To be sure, the fact that modern algebra became an autonomous discipline associated to the professional identity of being an "algebraist" during the inter-war period impacted retrospectively the history of "algebra." A key aspect of this new professional identity was the claim of the modernity of a conceptual approach to mathematics as emerging from a line of developments that involved especially Galois, Jordan, Dedekind, Hilbert, Noether, and Artin. For instance, Claude Chevalley argued that Artin was an algebraist because of his "intellectual temper" that consisted in focusing on "structures" and in "putting things in order" [Chevalley 1964]. When explaining why he considered Noether as "having been both the mother and the queen of the great German algebraic school," Paul Dubreil insisted on the conceptual nature of Noether's approach "in which computations are replaced by ideas" [Dubreil 1986, p. 16], thereby directly pointing to Galois's claim that one should jump over computations [Galois 1962, p. 9].

---

[74] See especially the way Julia opposed the model of the *grands hommes* to the cynic attitude and dandyism of the new generation as well as the way Cassou promoted the political involvement of the "citizen mathematician" Châtelet [Goldstein 2009, p. 173].

But we have seen that none of these claims were specific to the interwar modern algebra. A characterization of the algebraist-self very close to Chevalley's one had indeed already been given by Picard in 1890. Moreover, Dubreil's claim is almost the same as Adhémar's 1922 characterization of the "algebra of order" as quoted in the introduction of the present paper. As has been seen before, the "science of order" had been traditionally associated in France to Descartes's "modern algebra" as later developed by Poinsot who had made in 1808 a statement in regard to computations very similar to Galois's one [Boucard 2011].

It must therefore be pointed out that the categories used in public discourses on mathematics have some autonomy in regards with the evolutions of mathematics. In contrast, these categories depend of the targeted audiences. For instance, when he addressed a large audience with a topic as popular as E. Noether's biography, Dubreil claimed that modern algebra had been developed in Germany in the legacies of Galois and Jordan. When he dealt with more local topics, however, Dubreil not only pointed to some French algebraists of the beginning of the 20$^{th}$ century, such as de Séguier and Châtelet, but he also mentioned some mathematical themes he would otherwise have stayed silent about (such as Hermite's tradition in the theory of algebraic forms).

Picard's discourses provide an even more striking example. We have seen, that Picard had repeatedly advocated publicly the crucial role Jordan had played in the unfolding of the group theoretical ideas of Galois. But in his *Traité d'Analyse,* Picard nevertheless grounded his presentation of Galois theory on Kronecker's approach, which the latter had explicitly opposed to Jordan's one.

The complex interactions between public discourses and the evolution of mathematics raise the difficult issue of the impact of public discourses on the historiography of mathematics. Public discourses on Galois indeed seem to have had a lasting impact on the main categories of the historiography of algebra. Among these, one finds the implicit statement that national categories play a relevant role for characterizing disciplinary categories. But as has been alluded to in the first section of this papers, these two categories were already much interlaced in the late 19$^{th}$ century when histories of "disciplines" were developing at the same time as nation-state organisations.

Let us consider two additional examples. In the introduction of his textbook on group theory, Burnside insisted that after the foundational works of Galois and Cauchy, and despite the additions that were made by "French mathematicians" in the middle of the century, "no considerable progress in the theory, as apart from its applications, was made till the appearance in 1872 of Herr Sylow's memoir" [Burnside 1897, p. v]. But Jordan was both French and actually publishing his works on substitutions in the early 1870s. This example shows that the national category emphasized by Burnside had its own autonomy as a public category in regard with historical chronology. The public dimension of some categories sheds some light on the multiplicity of chronologies that often arise when one considers public discourses as testimonies, i.e., as historical sources [Goldstein 2009, p. 147].

The second example is the one of the national and cultural categories the reviewers of the *Jahrbuch* often emphasized in connection with disciplines or theories. For instance, Alfred Loewy claimed in 1901 that the issues tackled by Leonard Dickson's 1901 monograph on linear groups in Galois fields were mostly cultivated by English speaking people. George Miller swiftly reacted by publishing a paper in which he proved Loewy wrong by appealing to some French authorities such as Picard or Poincaré. Miller eventually concluded that "it need scarcely be added that some modern mathematicians seem to avoid group theory even where it would simplify the treatment of the subject in hand. This seems to be true, for instance, of Hilbert"s *Grundlagen der Geometrie.*" [Miller 1903, p.89]. Here Miller's claim was falling under Picard's authority. It is therefore in complete opposition to the usual focus of the historiography of algebra on Hilbert's Göttingen. But it is highly significant that both claims actually appeal to similar articulations between national/ disciplinary categories.

That this articulation exerted a long-term influence on the historiography of mathematic is further exemplified by the dichotomy between substitution groups and abstract groups. This dichotomy has indeed been a structuring one in the historiography of algebra in the 20$^{th}$ century.[75] But as an actor's category, it was much connected to national issues. For instance, as a *Jahrbuch* reviewer, Loewy contested the relevance of Séguier's 1904 monograph on abstract groups in claiming that, in contrast with the Germans, the French were not known for their contribution to contemporary abstract group theory but for some older works on substitution groups. But when one takes a closer look at Séguier's book, neither the categories of France nor algebra appear to have any obvious meaning. Indeed, on the one hand Séguier appealed to the works of the German Cantor and Frobenius. On the other hand the "abstract" approach Séguier promoted was in accordance with some contemporary works in analysis. Recall that even though the analysis of real function theory was much promoted in public discourses in France at the turn of the century, analysis was traditionally considered as including algebra and arithmetic through unifying concepts such as groups and functions. For instance, Hadamard celebrated in 1906 the "abstract" approach of Frechet's doctoral thesis [Goldstein 2009, p. 153], which the latter had explicitly connected to Séguier's work on abstract groups. In both Fréchet's abstract spaces and Séguier's abstract groups, the category "abstract" pointed to properties independent from the nature of the elements. This approach was thus quite opposed to the focus on the notion of Körper in the Göttingen-like algebraic number theory, but it was very close to the notion of Galois field of the Chicago algebraic school [Brechenmacher 2011].

Despite the evolution of the meanings taken on by national and disciplinary categories, the uses of such categories in public discourses show some long-term continuity. The impact of this situation on the historiography should be studied further. For instance, when Broglie presented Picard's achievements in 1942, the evolutions of analysis, arithmetic, and algebra as both autonomous specific disciplines and professional identities resulted in a radical change of the appreciation of Picard's involvement in "analysis." While the latter had claimed the unity of mathematics through analysis, he was now considered as having been "almost exclusively an analyst [...] Picard may have solved important problems in other branches of pure mathematics, such as Geometry and Arithmetic, but he always considered such problems through Analysis [Broglie 1942, p. 1].

---

[75] Cf. [Wussing, 1984]

## Conclusion

We have seen that most public discourses on Galois have underlying them a dichotomy between the world of pure ideas and the real world. This dichotomy points to the constitution of a public expression of mathematics in a tension between the legitimacies of the inner circle of the academy and the public opinions of some other practitioners of mathematics. Already in the early 1830s, the human rights for experiencing pity and respect had been opposed to the legitimacy of the experts. Later on at the turn of the century, Picard claimed the unifying power of Galois's groups for mathematics while Winter aimed at making Galois theory a part of the elite's culture. By then, the story of Galois had already given rise to a condensed - though extremely rich – tale of both universality and fragmentation. It points to the two-fold long-term process of autonomization of mathematics and of fragmentation of both the audiences and the practitioners of mathematics.

The issue of the fragmentation of the unity of mathematics was much discussed at the turn of century. Peano and Mittag-Leffler especially reflected upon the notions of "contributions," "being read and noticed," and how, while most mathematical issues "mean nothing to the general public," one could "construct and assemble a special public" in the international space [Turner 2011]. The present case study has appealed to a specific definition of the dichotomy between "public" and "mathematical" texts in connection with discourses on Galois. It calls for further investigations on the variety of identities taken on by mathematics for various audiences and in various time-periods.

Moreover, the tension between the legitimacies of academic authorities and of public opinions calls in turn for some further investigations on the transfers that occur between various fields of knowledge or activities through heterogeneous corpora. Indeed, several recent researches in the history of mathematics have shown the roles played by the intertextual relations of some collective organizations of texts, i.e., networks of texts, which identities are complex and often do not coincide with any disciplinary, national, or institutional identity. These investigations have changed our understanding of what are the relevant collective dimension for analyzing mathematical texts. They thus also call for new investigations on the transfers between various copora and various fields.

One may also wonder about the roles some actors who did not have access to global public expressions on mathematics may have played in the stratification of mathematics and of its audiences. For instance, the actors who designated themselves as "algebraists" in France were mostly in periphery of the main centers (i.e., the *Académie*, the Sorbonne, the *É.N.S.)*.

They regularly intervened on fragmented audiences: Séguier reviewed the 1897 edition of Galois's works in the jesuit journal *Études*, Adhémar published an abstracted version of Dupuy's biography in *La revue du mois*, Autonne incorporated Jordan's approach in a series of popular university textbooks, Robert de Montessus de Ballorre enriched the notice on Galois in the extended French translation he gave of Rouse-Ball's *History of mathematics* etc. These discourses were both influenced by the authorities of French mathematics and appealing to other sources. As a matter of fact, the algebraists often presented themselves as go betweens in claiming they aimed at mediating some mathematical works published in foreign countries (especially Germany) even though their works were usually not less original than many others. They also usually showed strong international ideals, thereby pointing to an international space that did neither coincide to most French authorities' national understanding of the international (such as Picard or Borel), nor to the institutional space of congresses or unions, or to the elitist "research" space emphasised by journals such as *Acta mathematica*.

Moreover, some mediations were in turn acknowledged as original works and thereby mediated themselves to new audiences. For instance, in contrast to Paul Tannery who lowered the historical value of the notes Montessus had added to Rouse Ball's history of mathematics in his review in the *Bulletin des sciences mathématiques,* the reviewer of the *Bulletin of the American Mathematical Society* emphasised the originality of Montessus's notes. In addition to these two bulletins which had been founded for developing links between mathematicians in reporting on recent publications, both the mathematical works and public interventions of peripheric actors such a Montessus or Séguier were sometimes echoed in journals such as the *Revue de métaphysique et de morale, la Revue philosophique, la Revue thomiste, la Revue de synthèse, la Revue générale des sciences pures et appliquées*, etc.

The concrete forms of the tranfers between various fields should be studied further. The publications of some reproductions, comments, translations, reformulations had been indeed consubstantial to the periodical form since its development in the 17th century [Peiffer and Vittu 2008]. For instance, we have seen that Pierpont had regretted that Dupuy's biography had not been added to the 1897 reprinting of Galois's works because the *Annales scientifiques de l'École normale* was not a journal "everyone" could access easily. But Pierpont's review was nevertheless itself a link in the chain of transfers that made Galois's biography "public." This chain involved some reprinting, such as of both Dupuy's biography and of Bertrand's comments on this biography [Ehrhardt 2007, p. 653]. The chain also involved reformulations, such as Adhémar's paper in *La revue du mois* through which Galois's story eventually went under Beaunier's pen in *Le Figaro*.

The active roles taken on by some actors in such transfers should also be studied further. For instance, Picard's 1924 discourse for the 50th anniversary celebration of the *Société mathématique de France* was published in both the *Bulletin de la S.M.F.* - a journal mostly devoted to the publication of original mathematical papers and in the *Revue générale des sciences pures et appliquées* (with a new title). Moreover, the part of the discourse that mentioned Galois was reproduced in Picard's toast at the celebrations of the 50th anniversary celebration of the French society of Physics, which again was published on several supports, including a monograph. It must be pointed out that these echoes were not limited to public discourses on mathematics. They also concerned more specialized publications. For instance, in the 1930s, George Bouligand developed a new interpretation of Galois's works through what he designated as "causal analysis." His views were first published in specialized journals, such as the *Comptes rendus de l'Académie des sciences de Paris*, but they were also reformulated for other audiences, such as in the *Revue générale des sciences pures et appliquées.* They were then reviewed and discussed in *Mathesis*, a journal devoted to the teaching of mathematics, as well in the periodical of history of sciences *Isis* or in the journal of literature *La revue de Paris.*

Such echoes raise the issue of the interplays between authority, reproductions, assimilations, and originality in publications on/of mathematics. The public dimension of the values of "originality," "creativity," or "research" seems indeed to have been often underestimated. We have seen in this paper that such values were instrumental to the ways some authorities both established and legitimated some boundaries between mathematics and the outside of mathematics as well as some stratifications and hierarchies within mathematics. It was thus compulsory for authors to legitimate their contributions in connection to such values, i.e., either by acknowledging their own role as mediators in the circulation of knowledge; or on the contrary by claiming the genuine novelty of their own individual contribution.

Again, this complex phenomenon of echoes both concerns mathematical papers and discourses on mathematics. For instance the notice on "Galois" in the encyclopedic dictionary *Larousse mensuel illustré* of 1929 [Augé, 1929a, p. 846] was extracted from d'Ocagne's popularization book *Hommes et choses de sciences.* The section on Galois in this book was in turn a reformulation of Picard's public discourses. Later on in the 1930s it was reformulated again in a collective volume of general history, as well as in an essay published in the journal of literature *La Revue Belge,* side by side with short stories by H. G. Wells, Paul Claudel, and Jean Rostand [Ocagne 1936, p. 436].[76]

Among the questions raised by echoes and transfers from one field to another, the question of the transfers of specific relations to History from a field to another is not the least challenging one. Indeed the value of originality in mathematics is closely connected with the promotion of dealng with free hands with referrences to previous works, i.e., with a kind of freedom in regard with History. A book, a paper, a proof, a theorem, or even a definition may be reformulated, cut to pieces, reorganized, and thereby disconnected from the issues it was originally devoted to. This process appears to be quite similar to the one of canibalism through which parts and parcels of old buildings or of outdated technical devices are both dismembered and being used for new purposes.

Even though selective references do not destroy physically the texts they target, they may nevertheless cause some parts of these texts to fall into oblivion. For instance, Jordan attributed to Galois a theorem the latter had never stated while the former did not even allude to the criterion the latter had concluded his main memoir with [Brechenmacher 2011]. While this criterion had been celebrated publicly for decades, it then fell into oblivion where it was soon followed by most other parts of Galois's works (especially modular equations and the analytic representations of substitutions). In contrast, some selected quotations have been well preserved. They thus keep standing in all their vainglory [Beaunier 1908], as an antique sphynx alone in a field of antique ruins.

---

[76] In this paper d'Ogagne actually appealed to genius scientists such as Galois to delimit the boundaries of the "feminine intelligence."

In claiming their creative power, mathematicians have also claimed the superiority of the mathematical truth over historical truth, i.e., their legitimacy for filling the holes in historical sources by revealing their thrue mathematical grounds. Authors, such as Jordan, Weber, Hilbert or Artin claimed they would reveal the "true" grounds of Galois theory while none of them made any direct reference to Galois's papers. It is therefore striking that the presentations of Galois's life in the fields of litterature, philosophy, history, politics, etc., have so often appealed to fiction to fill the holes in the historical sources on Galois's life [Albrecht & Weber 2011], thereby making a claim on the legitimacy of a fictional truth.

**Bibliography.**

ABEL (Niels Henrik),
[1826] Démonstration de l'impossibilité de la résolution algébrique des équations générales qui passent le quatrième degré, *Journal für die reine und angewandte Mathematik*, t. 1, 1826, p. 65-84.
[1829] Mémoire sur une classe particulière d'équations résolubles algébriquement, *Journal für die reine und angewandte Mathematik*, t. 4, 1829, p. 131-156.
 [1830a] Mathematische Bruchstücke aus Abels Briefen. *Journal für die reine und angewandte Mathematik* 5, 336-343.
[1830b] Fernere mathematische Bruchstücke aus Abels Briefen, Schreiben Abels an Legendre. *Journal für die reine und angewandte Mathematik* 6, 73-80.
ADHEMAR (Robert d'),
[1905] Trois maîtres : Ampère, Cauchy, Hermite, *La Quinzaine,* n°265, 1905/11/01-1905/12/16, p. 1-16.
[1922] Nécrologie. Camille Jordan, *Revue générale des sciences pures et appliquées*, t. 3 (15 février), 1922, p. 65-66.
[1924] Election de M. Picard à l'Académie française, *Revue générale des sciences pures et appliquées*, n°23 (1924), p. 623-624.
ALAIN,
[1909] Évariste Galois, *La dépêche de Rouen*, 10 août 1909
ALBRECHT (Andrea) et WEBER (Anne-Gaëlle),
[2011] La réception littéraire de Galois, *Revue d'histoire des mathématiques*, to appear.
L'ALBUM NATIONAL (Anonyme),
[1829a] [Recension sur] Académie des sciences. Séance du lundi 2 mars, *L'album national, revue de la littérature, des sciences, des cours publics, des tribunaux, des théatres, des arts et des modes,* 7 mars 1829, p. 260.
[1829b] [Recension sur] Académie des sciences. Séance du lundi 22 juin, *L'album national, revue de la littérature, des sciences, des cours publics, des tribunaux, des théatres, des arts et des modes,* 27 juin 1829, p. 516.
ANDRE (Aimé), éd.
[1834] Evariste Galois, *Dictionnaire biographique universel et pittoresque, contenant 3000 articles environ de plus que la plus complète des biographies publiées jusqu'a ce jour,* Paris : Aimé André, 1834, p. 464*.*
L'ANNEE SCIENTIFIQUE ET INDUSTRIELLE (Anonyme)
[1883] Nécrologie scientifique. Liouville, *L'Année scientifique et industrielle.* 1883, p. 524-526.
APOLOGETIQUE (Anonyme)
[1915] La science et la philosophie allemande, *Revue pratique d'apologétique,* 1915/02, p. 365-372.
ARCHIBALD (Thomas)
[2011] Differential equations and algebraic transcendents: French efforts at the creation of a Galois theory of differential equations (1880-1910), *Revue d'histoire des mathématiques,* à paraître.
ARNAULT (Antoine Vincent), BAZOT (Étienne-François), et al.,

[1825], *Biographie nouvelle des contemporains ou Dictionnaire historique et raisonné de tous les hommes qui, depuis la Révolution française, ont acquis de la célébrité par leurs actions, leurs écrits, leurs erreurs ou leurs crimes, soit en France, soit dans les pays étrangers*, Paris : Librairie historique, 1820-1825.

ARTIN (Emil)
[1938] *Foundations of Galois Theory*, New York: New York University Lecture Notes, 1938.

AUBIN (David)
[forthcoming] War Cultures in Postwar French Mathematics: A Critique of the Bourbaki and Forman Theses, preprint.

AUBIN (David) and BIGG (Charlotte),
[2003] Neither Genius nor Context Incarnate : Norman Lockyer, Jules Janssen and the Astrophysical Self, *in* Thomas Söderqvist (éd.), *The History and Poetics of Scientific Biiography*, Londres, Ashgate Publilshing : 51-70.

AUBRY (A.),
[1894-1897] Essai historique sur la théorie des équations, J. de Math. spéc. (4) III (1894), p. 225-228, 245-253, 276-278, V (1897), p. 21, 17-20, 61-62, 83-88, 114-115, 131-132, 155-159.

AUDIGANNE (P.), BAILLY (EUGENE), CARISSAN (EUGENE), et al.
[1861] *Paris dans sa splendeur,* Paris : Charpentier, 1861, part. 1, volume 2.

AUGE (Paul) (dir.),
[1929a] [Résumé de] *Hommes et choses de sciences* de Maurice d'Ocagne. *Larousse mensuel illustré, revue encyclopédique universelle*, 1929-1931, p.844-846.
[1929b] Mathématiques. *Larousse mensuel illustré, revue encyclopédique universelle*, 1929-1931, p.624-626.

AULARD (F.A.),
[1897] Chronique et bibliographie, *La révolution française. Revue d'histoire moderne et contemporaine,* t. 32 (1897), p. 86-87.

AUTONNE (Léon)
[1913] *Review on* De Séguier. Théorie des groupes finis, *Revue générale des sciences pures et appliquées,* t. 24 (1913), p. 401.

BEAULIEU (Lilianne)
[1993] A Parisian café and ten proto-Bourbaki meetings (1934- 1935), *The Mathematical intelligencer*, 15/1 (1993), 27-35.
[2009] Regards sur les mathématiques en France entre les deux guerres, *Revue d'histoire des sciences*, 62-1 (2009), p. 9-38

BEAUNIER (André)
[1908] A travers les revues. Évariste Galois. *Le Figaro littéraire.* October, 17, 1908, p. 3.

BELHOSTE (Bruno).
2001] The Ecole Polytechnique and Mathematics in Nineteenth-Century France, *in*
*Changing Images of Mathematics. From the French Revolution to the New Millenium*, Umberto
Bottazzini and Amy Dahan (ed.), London: Routledge, p. 15–30.

BELL (Eric Temple),
[1937] Genius and Stupidity. Galois, *Men of Mathematics*, New York, Simon and Schuster, 1937 p. 362-378.

BENSAUDE-VINCENT (Bernadette),
[2003] *La science contre l'opinion. Histoire d'un divorce,* Paris, Les Empêcheurs de penser en rond, le Seuil.

BENSAUDE-VINCENT (Bernadette) and RASMUSSEN (Anne) (dir.)[
[1997] *La science populaire dans la presse et l'édition, XIXe et XXe siècles,* Paris, CNRS Éditions : 51-68.

BERTRAND (Joseph),
[1867] *Rapport sur les progrès les plus récents de l'analyse mathématique*, Ministère de l'instruction publique, Imprimerie impériale, Paris, 1867.
[1899] La vie d'Évariste Galois par P. Dupuy, *Journal des savants*, juillet 1899, p. 389-400 [repr. *Éloges Académiques*, Paris : Hachette, 1902, p.329-345].

BIBAS (Henriette),
[1947] Deux notes sur le journal de Vigny et sur "Volupté", *Société d'histoire littéraire de la France. Revue d'histoire littéraire de la France*, 1947/01/01-1947/03/31, p. 271-274.

BLANC (Louis),
[1840] La Commune. *Revue du progrès politique, social et littéraire*, vol. 4 à 5, p. 46-66.
[1842] *Histoire de dix ans*, Paris, Pagnerre, 1842.

BIBAS (Henriette)
[1947] Deux notes sur le journal de Vigny et sur Volupté, Société d'histoire littéraire de la France. Revue d'histoire littéraire de la France. 1947/01/01-1947/03/31, p. 272-274.

BIRKHOFF (Garret)
[1937] Galois and Group Theory, *Osiris*, vol. 3, 1937, p. 260-268.
BÔCHER (Maxim)
[1914] Charles Sturm et les mathématiques modernes, *Revue du mois,* 17 (1914), p. 88-104.
BOLTANSKI (Luc) and THÉVENOT (Laurent)
[1991] *De la justification. Les économies de la grandeur,* Paris, Gallimard, 1991.
BONNET (Jean-Claude)
[1986] Les morts illustres : oraison funèbre, éloge académique, nécrologie, in Pierre Nora (éd.), *Les lieux de mémoire, t. II, La nation,* vol. 3 Paris, Gallimard : 217-241
BOREL (Émile),
[1904] Les exercices pratiques de mathématiques dans l'enseignement secondaire, Revue générale des sciences pures et appliquées, 14 (1904), p. 431-440.
BOREL (Émile), DRACH (Jules)
[1995] *Introduction à l'étude de la théorie des nombres et de l'algèbre supérieure*, Paris : Nony, 1895.
BOYER (Jacques),
[1900] *Histoire des mathématiques*, Paris : Gauthier-Villars, 1900.
BOUCARD (Jenny),
[2010] Louis Poinsot et la théorie de l'ordre : un chaînon manquant entre Gauss et Galois ?, *Prépublication.*
BOULIGAND (Georges),
[1934] Remarques sur l'analyse causale des théorèmes de géométrie, *Revue générale des sciences pures et appliquées,* t. 45, 1934, p. 322-323.
BOURBAKI (Nicolas),
[1960] *Eléments d'histoire des mathématiques,* Paris : Hermann, 1960.
BOURDIEU (Pierre)
[1976] Quelques propriétés des champs, *in Questions de sociologie, Exposé à l'E.N.S., novembre 1976*, Paris 2002 : Les Editions de minuit, p. 113-120
[2001] *Langage et pouvoir symbolique*, Paris, Seuil, 2001.
BOURGNE (Robert), AZRA (Jean-Pierre),
[1962] *Écrits et mémoires mathématiques ; édition critique intégrale des manuscrits et publications d'Évariste Galois par Robert Bourgne et Jean-Pierre Azra*, Paris, Gauthier-Villars, 1962.
BOURLET (Carlo),
[1907] L'enseignement de la géométrie *Bulletin de la Société française de Philosophie*, séance du 21 mars 1907.
BOURQUELOT (Félix), MAURY (Alfred),
[1852] Evariste Galois, *La littérature française contemporaine, 1827-1849, continuation de la France littéraire. Dictionnaire bibliographique,* vol. 4., 1852, p. 17.
[1854] Joseph Liouville, *La littérature française contemporaine, 1827-1849, continuation de la France littéraire. Dictionnaire bibliographique,* vol. 5., 1854, p. 161.
BOUTROUX (Pierre),
[1913] [Recension sur] Léon Brunschwig, *Les étapes de la philosophie mathématique*, Revue de métaphysique et de morale, 1913, p. 107-131.
[1919] L'histoire des sciences et les grands courants de la pensée mathématiques*, La revue du mois,* t.20 (1919, p. 604-621.
[1920] *L'idéal scientifique des mathématiciens dans l'antiquité et les temps modernes,* Paris : Alcan, 1920.
BOYER (Jacques)
[1900] *Histoire des mathématiques*, Paris : Carré et Naud, 1900.
BRECHENMACHER (Frédéric).
[2007] La controverse de 1874 entre Camille Jordan et Leopold Kronecker, *Revue d'Histoire des Mathématiques*, tome 13, fasc. 2, p. 187-257.
[2011] Self-portraits with Evariste Galois (and the shadow of Camille Jordan), *Revue d'histoire des mathématiques*, to appear.
[2012] On Jordan's measurements*, to appear.*
BROGLIE (Louis de),
[1942] *La vie et l'oeuvre d'Émile Picard*, Paris, Institut de France, Gauthier-Villars, 1942.
BROKS (Peter)
[1996] *Media Science before the Great War*, London, MacMillan Press, 1996.
BRUNSCHWIG (Léon),
[1912] *Les étapes de la philosophie mathématique*, Paris, Alcan 1912.

[1923] [Recension sur] Pierre Boutroux, *L'idéal scientifique des mathématiciens dans l'antiquité et les temps modernes*, *Revue philosophique de la France et de l'étranger. 1923 . Janv.-juin,* p. 154-158.
    BURNSIDE (William)
[1897] *Theory of Groups of Finite Order*, Cambridge: Cambridge University Press, 1897.
    CAHEN (Eugène)
[1900] *Éléments de la théorie des nombres*, Paris: Gauthier-Villars, 1900.
    CAJORI (Florian),
[1919] *A History of Mathematics*, 2e éd., New York : MacMillan, 1919.
    CARNOT (Hyppolite) et LEROUX (Pierre),
[1832] Travaux mathématiques d'Evariste Galois, *Revue encyclopédique*, 1932, p. 566-568.
    CHARLE (Christophe)
[1994] *La république des universitaires 1870-1945,* Paris, Seuil, 1994.
    CHEVALLEY (Claude),
[1964] Emil Artin (1898-1962), *Bulletin de la société mathématique de France,* 92, p. 1-10.
    CHEVALIER (Auguste),
[1832] Nécrologie : Evariste Galois, *Revue encyclopédique*, t.55, juill.-sept. 1832, p. 744-754.
    CHEVREUSE (Louis)
[1912] Le manuscrit, *Le figaro,* 10/4/1912, p.1.
    COMBEROUSSE (Charles de)
[1898] *Cours d'algèbre supérieure,* Deuxième édition, Paris : Gauthier-Villars, 1898.
    CORRY (Leo)
[1996] *Modern Algebra and the Rise of Mathematical Structures*, Basel : Birkhäuser, 1996.
    COUTURAT (Louis)
[1898] Sur les rapports du nombre et de la grandeur, *Revue de métaphysique et de morale*, VI (1898), p. 422-447.
    CURIE (Marie)
[1924] *Pierre Curie,* Paris : Payot, 1924.
    DARBOUX (Gaston)
[1911] *Éloge des donateurs de l'Académie*, Paris, Institut de France, Gathier-Villars, 1911.
    DASTON (Lorraine), SIBUM (Otto)
[2003] Introduction : Scientific Personae and Their Histories, *Science in Context*, vol. 10, n°2 : 95-113.
    DENJOY (Arnaud)
[1934] Paul Painlevé, *Annales de l'Université de Paris*, 1934, p. 1-23.
    D'ESCLAYBES (Robert),
[1898] Revue des livres : *Traité d'algèbre supérieure. Compagnie de Jésus. Études de théologie, de philosophie et d'histoire*, 1898, p. 414-416.
    DESCHANEL (Paul),
[1919] Les Allemands et la science, in *La France victorieuse, paroles de guerre*, Paris : 2e mille, 1919, p. 118-131.
    DESROCHERS (Pierre-Charles),
[1833] Evariste Galois, *Nécrologie de 1832 ou notices historiques sur les hommes les plus marquans tant en France que dans l'étranger, morts pendant l'année 1832,* Paris, 1932, p. 132.
    DIEUDONNE (Jean),
[1962] Notes sur les travaux de Camille Jordan relatifs à l'algèbre linéaire et multilinéaire et la théorie des nombres, [Jordan, *Œuvres*, 3, p. V-XX].
    DUBREIL (Paul).
[1982] L'algèbre, en France, de 1900 à 1935, *Cahiers du séminaire d'histoire des mathématiques*, 3 (1982), p. 69-81.
[1986] Emmy Noether, *Cahiers du séminaire d'histoire des mathématiques*, 7 (1986), p. 15-27.
    DUCLERT (Vincent), RASMUSSEN (Anne)
[2002] La république des savants, *in* Duclert (Viencent), Prochasson (Christophe), dir., *Dictionnaire critique de la République*, Paris, Flammarion : 439-445.
    DUFUMIER (H.),
[1911] La généralisation mathématique, *Revue de métaphysique et de morale*, 1911, p. 723-752.
    DUNAND (Renée),
[1923] Les derniers livres parus, *La Pensée française. Organe d'expansion française et de propagation nationale*, 1923/12/13, p. 17-18.
    DUPUY (Paul),

[1896] La vie d'Évariste Galois, *Annales scientifiques de l'École normale supérieure*, 3e sér., t. 13, 1896, p. 197-266.
    EHRHARDT (Caroline).
[2007] *Evariste Galois et la théorie des groupes. Fortune et réélaborations (1811-1910)*, Thèse de doctorat. Ecole des Hautes études en sciences sociales. Paris, 2006.
[2010a] A Social History of the Galois's Affair at the Paris Academy of Sciences, *Science in Context*, 23(1), 2010, p. 91-119.
[2010b] La naissance posthume d'Évariste Galois, *Revue de synthèse*, *t.* 131, 6e série, n° 4, 2010, p. 543-568
    ELIAS (Norbert)
[1982] "Scientific establishments," *in* id., Herminio Martins, Richard Whitley (éd.), *Scientific Establishments and Hierarchies*, Dordrecht, Reidel, 1982 : 3-69.
D'ESCLAYBES (Robert)
    [1898], Revue des livres : *Traité d'algèbre supérieure. Compagnie de Jésus. Études de théologie, de philosophie et d'histoire*, 1898, p. 414-416.
L'EUROPE NOUVELLE (Anonyme),
[1918] Les sciences, *L'Europe nouvelle, revue hebdomadaire des questions extérieures, économiques et littéraires*, n°23, 15 juin 1918, p. 1109-1110.
    FEHR (H.),
[1897] [Recension sur] Petersen (Julius), Théorie des équations algébriques*, Revue générale des sciences pures et appliquées*, t. 8 (1897), p.756.
    FELLER (François Xavier de)
[1832-1836] Evariste Galois, *La Biographie universelle, ou Dictionnaire historique des hommes qui se sont fait un nom par leur génie, leurs talents, leurs vertus, leurs erreurs ou leurs crimes*, Lille, 1932, vol 13 p. 93, réed. 1834, vol5, p. 533, réed. 1836, vol. 9, p. 53.
    LE FIGARO (Anonymous),
[1924] Les trois élections d'hier à l'Académie, *Le Figaro,* 28/11/1924, p.1.
[1941] Un grand philosophe mathématicien disparaît. *Le Figaro,* 16/12/1941, p.4.
    FRANCE ET MONDE (Anonymous),
[1922] La documentation vivante. L'histoire des sciences et les prétentions de la science allemande, par Emile Picard, *France et monde. Revue de documentation économique et sociale*, 1922/09/20., p. 304.
    CASTELLI GATTINARE (Enrico)
[2001] Épistémologie 1900. La tradition française, *Revue de synthèse*, 4e sér., n° 2-3-4, 2001, p. 347-365
    GISQUET (Henri)
 [1840] *Mémoire de ma vie*, Paris: Marchant, 1840.
    GALOIS (Evariste),
[1846] Œuvres mathématiques*, Journal de mathématiques pures et appliquées* 11 (1846), 381–444.
[1897] *OEuvres mathématiques d'Évariste Galois, publiées sous les auspices de la Société Mathématique de France, avec une introduction par M. Émile Picard*, Paris Paris: Gauthier-Villars, 1897.
[1906] Manuscrits et papiers inédits de Galois ; par M. Jules Tannery, 1re partie », *Bulletin des sciences mathématiques*, 2e sér., t. 30, , 1906, p. 226-248 et p. 255-263.
[1907] Manuscrits et papiers inédits de Galois ; par M. Jules Tannery, 2e partie », *Bulletin des sciences mathématiques*, 2e sér., t. 31, 1907, p. 275-308.
 [1962] *Ecrits et mémoires mathématiques*, ed. R. Bourgne, J.-P. Azra. Paris: Gauthier-Villars, 1962.
    GAUTHIER (Sébastien)
 [2009] La géométrie dans la géométrie des nombres : histoire de discipline ou histoire de pratiques à partir des exemples de Minkowski, Mordell et Davenport, *Revue d'histoire des mathématiques,* **15-2** (2009), p. 183-230.
    GISPERT (Hélène),
[1991] La France mathématique*. La Société Mathématique de France (1870-1914)*, Cahiers d'histoire et de philosophie des sciences, Paris : Belin, 1991.
[1995] La théorie des ensembles en France avant la crise de 1905 : Baire, Borel, Lebesgue... et tous les autres, *Revue d'histoire des mathématiques*, I (1995), p. 39-81.
    GISPERT (Hélène) and TOBIES (Renate)
[1996] A comparative study on the French and German Mathematical Societies before 1914, [Goldstein, Gray, Ritter 1996].
    GOLDSTEIN (Catherine)
 [1999] Sur la question des méthodes quantitatives en histoire des mathématiques : le cas de la théorie des nombres en France (1870- 1914) », *Acta historiae rerum necnon technicarum*, nouv. sér., vol. 3, 1999, p. 187-214.

[2009] La théorie des nombres en France dans l'entre-deux-guerres : De quelques effets de la première guerre mondiale, *Revue d'histoire des sciences* 62(1), 2009, 143-176.
[2011] Hermite's strolls in Galois fields, *Revue d'histoire des mathématiques Revue d'histoire des mathématiques*, t. 17, fasc. 2, *to appear*.
  GOLDSTEIN (Catherine), GRAY, (Jeremy), RITTER (Jim) (eds.)
[1996] *L'Europe mathématique. Histoires, mythes, identités*, Paris : Éditions de la MSH, 1996.
  GOLDSTEIN (Catherine), SCHAPPACHER (Norbert), SCHWERMER (Joaquim)(eds.),
[2007] *The Shaping of Arithmetics after C. F. Gauss's Disquisitiones Arithmeticae*, Berlin : Springer, 2007.
  GOMPERZ (Theodor)
[1908] *Les penseurs de la Grèce : histoire de la philosophie antique*, vol III, Paris : Alcan, 1910.
  GUÉRIN (Charles)
[2007] L'élaboration de la notion rhétorique de le *persona* au 1er siècle avant J.C.: antécédents grecs et enjeux cicéroniens, *L'information littéraire,* vol. 59 (2007), n°2 : 37-42.
  GUYOT (Yves)
[1867] *L'inventeur*, Paris : Le Chevalier, 1867.
  HADAMARD (Jacques),
[1912] Le calcul fonctionnel, *L'enseignement mathématique*, 14 (1912), p. 5-19.
[1913a] L'œuvre d'Henri Poincaré. Le mathématicien, *Revue de métaphysique et de morale*, 1913, p. 617-658.
[1913b] Henri Poincaré et le problème des trois corps, *La Revue du mois*, 1913, p. 386-418.
[1920] Rapport sur les travaux examinés et retenus par la Commission de Balisitique pendant la durée de la guerre, *Comptes rendus de l'Académie des sciences de Paris*, t. 171 (1920), p. 436-445.
[1924] Les grandes idées humaines. La pensée française dans l'évolution des sciences exactes, *France et monde. Revue de documentation économique et sociale*, 1922/09/20, p. 321-343.
  HANNEQUIN (Arthur)
[1908] *Études d'histoire des sciences et d'histoire de la philosophie*, Paris : Alcan, 1908.
  HOVELAQUE (Émile)
[1911] Jules Tannery, *La revue de Paris*, 18-1 (1911), p.305-322.
  IHL (Olivier)
[2007] *Le Mérite et la République. Essai sur la société des émules*, Paris, Gallimard.
  L'HUMANITE (anonyme)
[1916] [Recension sur] *La science française,* L'Humanité. 23/04/1916, p. 3.
[1923] Les normaliens et la révolution. Evariste Galois. *L'humanité*, 02/06/1923, p. 6.
[1924] Les nouveaux élus de la vieille académie *L'humanité*, 28/11/1924, p. 2.
  JORDAN (Camille)
[1870] *Traité des substitutions et des équations algébriques,* Paris, 1870.
  JULIA (Gaston)
[1919] Lambert (Paul-Jean-Étienne), né à Annecy le 27 février 1894, tué à l'ennemi près de Fontenoy (Aisne) le 15 mars 1915. – Promotion de 1911, *in* Association amicale de secours des anciens élèves de l'Ecole normale supérieure (Paris), *Réunion générale annuelle* (1919), 109-113.
  KIERNAN (Melvin),
[1971] The Development of Galois Theory from Lagrange to Artin, *Archive for History of Exact Sciences*, vol. 8, n° 1-2, 1971, p. 40-152
  LACROIX (Sylvestre-François) and POISSON (Siméon Denis)
[1831] Rapport sur le mémoire de M. Galois relatif aux conditions de résolubilité par radicaux, *Procès-verbaux des séances de l'Académie* (4 juillet 1831), p. 660-661.
  LAFITTE (Jean-Paul)
[1916] La science française, *L'humanité,* 29 avril 1916, p. 3.
  LAROUSSE (Pierre) dir.,
[1866-1877a] Galois, *Grand Dictionnaire universel du XIXe siècle, français, historique, géographique, mythologique, bibliographique, littéraire, artistique, scientifique, etc., etc.,* Paris, 1872, p. 549.
[1866-1877b] Liouville, *Grand Dictionnaire universel du XIXe siècle, français, historique, géographique, mythologique, bibliographique, littéraire, artistique, scientifique, etc., etc.,* Paris, 1873, p. 549.
  LASSEUR,
[1933] Influence de l'âge sur la personnalité scientifique. Le laboratoire du professeur honoraire. *Revue générale des sciences pures et appliquées*, 1933, p. 647-649.
  LEBESGUE (Victor Amédée),
[1846] Sur l'inscription des polygones réguliers de 15 et de 17 côtés, *Nouvelles annales de mathématiques,* 1re série, t. 5 (1846), p. 683-689.

   LEFEVRE, (B.),
[1897] *Cours développé d'algèbre élémentaire précédé d'un aperçu historique sur les origines des mathématiques élémentaires et suivi d'un recueil d'exercices et de problèmes*. Wesmaer-Charlier, Namur, 1897.
   LELOUP (Juliette),
[2009] *L'entre-deux-guerres mathématique à travers les thèses soutenues en France*, thèse de doctorat, Université Pierre et Marie Curie, 2009.
   LIE (Sophus),
[1895] Influence de Galois sur le développement des mathématiques, in DUPUY (Paul) (ed.), *Le Centenaire de l'École Normale 1795-1895*, Paris, Hachette, 1895.
   LIEBER (Lillian),
[1932] *Galois and the Theory of Groups : a Brigtht Star in Mathesis*, Lancaster, The Science Press Printing Company, 1932.
   LIOUVILLE (Joseph),
[1846] « Avertissement », in [Galois, 1846, p. 382-384].
   LOVE (A.E.H.)
[1917] La recherche mathématique, *Revue générale des sciences pures et appliquées,* 28 (1917), p. 271-275
   MAGASIN PITTORESQUE (Anonyme)
[1848] Évariste Galois, *Magasin pittoresque*, t. 16, 1848, p.227-228.
   MANSION (Paul)
[1910] La légende de Galois, *Annales de la société scientifique de Bruxelles*, 1910, p. 104-105.
   MARIANI (J.)
[1932] Evariste Galois et l'évolution des mathématiques. *Centre international de Synthèse, Berr, Henri (dir.). Revue de Synthèse,* 1932, p. 7-14*.*
   MAZLIAK (Laurent)
[forthcoming] The ghosts of the École normale. Life, death, and destiny of René Gateaux, *preprint.*
   MAUSS, Marcel
[1938] Une catégorie de l'esprit humain : la notion de personne, celle de "Moi," *The journal of the Royal Anthropologica Institute ofGreat Britain and Ireland,* vol 68 (1938): 263-281.
   MEHRTENS (Herbert)
[1990] *Moderne - Sprache - Mathematik. Eine Geschichte des Streits um die Grundlagen der Disziplin und des Subjekts formaler Systeme*. Frankfurt : Suhrkamp,, 1990.
   MENTRE (F.),
[1919] Les lois de la production intellectuelle. *Revue philosophique de la France et de l'étranger. 1919 . Juil.-déc.p.* 447-478.
   MILLER (George Abram),
[1903] Appreciative Remarks on the Theory of Groups, *The Amer. Math. Monthly*, vol. 10, 1903, p. 87-89
   MOUREU (Charles),
[1928] La science française, *Revue politique et parlementaire*, t. 134, p. 331-349.
   NABONNAND (Philippe),
[2000] La polémique entre Poincaré et Russell au sujet du statut des axiomes de la géométrie, *Revue d'histoire des mathématiques*, t.6, 2 (2000), p. 219-269.
   NERVAL (Gérard de)
[1841] Mémoire d'un parisien, *L'artiste*, 11 avril 1841
   NICOLET (Claude)
[1994] *L'idée républicaine en France (1789-1924). Essai d'histoire critique,* Paris, Gallimard, 1994
   NYE (Mary Jo)
[1979] The Boutroux Circle and Poincaré's Conventionalism, *Journal of the History of Ideas*, vol. 40, n° 1, 1979, p. 107-120
   D'OCAGNE (Maurice),
[1930-1936] *Hommes* & *choses de science - Propos familiers*,  Paris : Librairie Vuibert, 3 Vol., 1930-1936.
[1934] Un coin de la France Intellectuelle, *La revue de Paris*, sept. oct., 1934, p. 391-410.
[1936] L'intelligence féminine, *La revue belge*, p. 433-443.
   PANNIER (Sophie)
[1831-1834] Un jeune républicain en 1832 », in *Paris, ou le livre des Cent-et-un*, Paris, Ladvocat, 1831-1834, t. 10, p. 197-217.
   PEIFFER (Jeanne), VITTU (Jean-Pierre),
[2008] Les journaux savants, formes de la communication et agents de la construction des savoirs (17é-18e siècles), *Dix-huitième siècle*, 40 (2008), p. 241-259.

   PHILIPPE (Léon),
[1913] Le progrès musical, *Annales de l'Institut international de sociologie*, 1913, p. 325-344.
   PICARD (Émile),
[1883] Sur les groupes de transformation des équations différentielles linéaires », *Comptes rendus hebdomadaires des séances de l'Académie des sciences*, vol. 46, 1883, p. 1131-1134.
[1890a] Revue annuelle d'analyse, *Revue générale des sciences pures et appliquées*, 1 (1890), p. 702-708.
[1890b] Notice sur la vie et les travaux de Georges-Henri Halphen, *Comptes rendus de l'Académie des sciences de Paris,* t. 111 (1890), p. 489-497.
[1896] *Traité d'analyse*, t. III, Paris, Gauthier-Villars, 1896.
[1900] L'idée de fonction depuis un siècle, *Revue générale des sciences pures et appliquées*, 1900, p. 61-68.
[1902] Mathématiques, *in* Alfred Picard (éd.), *Exposition universelle internationale de 1900. Le bilan d'un siècle (1801-1900),* Imprimerie nationale : Paris, t.1, 1902, p. 121-137.
[1914] *La science moderne et son état actuel*, Flammarion : Paris, 1914.
[1916] *L'histoire des sciences et les prétentions de la science allemande*, Perrin et cie : Paris, 1916.
[1922a] Résumé des travaux mathématiques de Jordan, *Comptes rendus de l'Académie des sciences de Paris,* t. 174 (1922), p. 210-211.
[1922b] *Discours et mélanges*, Gauthier-Villards : Paris, 1922.
[1924a], Discours de M. Emile Picard, *Bulletin de la S.M.F.,* t. 52, 1924, p.27-32.
[1924b] De l'objet des sciences mathématiques, *Revue générale des sciences pures et appliquées*, 1924, p.325-327.
[1925] Toast de M. Emile Picard, *in Livre du Cinquantenaire de la Société Française de Physique*, Paris, 1925, p. 22-24.
   PIERPONT (James),
[1895a] Zur Geschichte der Gleichung des V. Grades (bis 1858), *Monatsh. f. Math.* VI (1895), p. 15-68.
[1895b] « Lagrange's Place in the Theory of Substitutions », *Bulletin of the American Mathematical Society*, vol. 1, 1894, p. 196-204.
[1897] Early History of Galois Theory of Equations, *Bulletin of the American Mathematical Society*, vol. 2, n° 4, 1897, p. 332-340.
[1899] [Review on] Galois' Collected Works. *Bulletin of the American Mathematical Society*. Vol. 5, N°6 (1899), p. 296-300.
   POINCARE (HENRI)
[1905] *La science moderne et son état actuel*, Flammarion : Paris, 1914.
[1910] *Savants et Écrivains*, Paris, Flammarion, 1910.
   LE PRÉCURSEUR (Anonymous)
[1832] Paris, 1er juin ; Correspondance particulière, *Le Précurseur. Journal constitutionnel de Lyon et du Midi*, 4-5 juin 1832.
   RASPAIL (François-Vincent),
[1829a] [Review on] Académie de sciences de Paris. 26 janvier. *Annale des sciences de l'observation* … par MM. Saygey et Raspail, T. 1, p. 336.
[1829b] Note du rédacteur. *Annale des sciences de l'observation* … par MM. Saygey et Raspail, T. 3, p. 144.
[1830] Côteries scientifiques. *Annale des sciences de l'observation* … par MM. Saygey et Raspail, T. 3, p. 150-152.
[1839] *Réforme pénitentiaire. Lettres sur les prisons de Paris*, Paris, Tamisey et Champion, 1839.
   REVUE DE METAPHYSIQUE ET DE MORALE (Anonymous),
[1921] [Recension sur] Pierre Boutroux, *L'idéal scientifique des mathématiciens*, *Revue de métaphysique et de morale*, 1921.p. 6-7.
   REVUE GENERALE DES SCIENCES PURES ET APPLIQUEES (Anonymous),
[1900] La bosse des mathématiques, *Revue générale des sciences pures et appliquées,* n°15 (1900), p.913-914.
   RIC ET RAC (Anonymous)
[1932] Un duel en 1932, *Ric et Rac : Grand hebdomadaire pour tous,* 1932/05/28, p. 2.
   RICHET (Charles),
[1923] *Le savant*, Paris : Hachette, 1923.
   ROUGIER (L.),
[1913] [Recension sur] Buche (Joseph), *La reprise de la querelle des anciens et des modernes*, *La phalange,* 1913, p. 321-339.
   ROUSE BALL (Willia Walter)
[1888] *A short account of the history of mathematics,* New York : Dover, 1888.
   RUPKE (Nicolaas A.),

[2008] *Alexander von Humboldt: A Metabiography*, U of Chicago P, 2008.
    SAIGEY (Jacques),
[1829] Abel, nécrologie, *Annales des sciences d'observation,* t. 2, 1829, p. 317-321.
    SARTON (George),
[1912] La chronologie de l'histoire de la science, *Revue générale des sciences pures et appliquées*, 1912, p. 341-342.
[1921] « Evariste Galois », *The Scientific Monthly*, oct. 1921, p. 363-375 [repr. in *Osiris*, vol. 3, 1937, p. 241-259].
    SÉGUIER (Jean-Armand de),
[1897] [Review on] *Œuvres mathématiques d'Evariste Galois*, Revue des livres, *Compagnie de Jésus. Études de théologie, de philosophie et d'histoire*, 1897, p. 139-140.
[1904] *Éléments de la théorie des groupes abstraits*, Paris : Gauthier-Villars, 1904.
    SERRET (Joseph Alfred),
[1849] *Cours d'algèbre supérieure*, Paris : Bachelier, 1849.
[1866] *Cours d'algèbre supérieure*, 3e éd. Paris, Gauthier-Villars, 1866, 2 vols.
    SYLOW (Ludwig) ,
[1920] Évariste Galois, *Norsk matematisk tidsskrift*, vol. 2, 1920, p. 1-17.
    TANNERY (Jules),
[1909] La vie et l'oeuvre d'Évariste Galois, *Revue scientifique*, t. 31 juillet 1909, p. 129-132
    TANNERY (Paul),
[1898] Les sciences en Europe, in *Histoire générale du IVe siècle à nos jours*, Lavisse (Ernest) et Rambaud (Alfred) (dir.), 1898, p. 731-767.
    TATON (René),
[1947] Les relations scientifiques d'Évariste Galois avec les mathématiciens de son temps, *Revue d'histoire des sciences*, t. 1, 1947, p. 114-130.
[1971] Sur les relations mathématiques d'Augustin Cauchy et Évariste Galois, *Revue d'histoire des sciences*, t. 24, 1971, p. 123-148.
[1993] Évariste Galois et ses biographes. De l'histoire aux légendes, *Sciences et techniques en perspective*, t. 26, 1993, p. 155-172.
    THIBAUDET (Albert)
[1930] Stendhal, le centenaire du *Rouge et Noir*, *La revue de Paris,* 11-12/1930, p. 317-336.
    TERQUEM (Olry),
[1849] Biographie. Richard, Professeur, *Nouvelles annales de mathématiques*, n° 3, 1849, p. 448-452.
    TRELAT
[1840] [Recension sur] *Des établissements d'éduction* de M. de Fellenberg, A. Bofwyl, traduction libre de l'allemand par M. Eugène de Caffarelli, maître de requêtes, *Revue du progrès politique, social et littéraire*, vol. 4 à 5, p. 110-119.
    TURNER (Laura E.)
[2011] *Identities, agendas, and mathematics in an international space*, PHD Thesis, Aarhus University, 2011.
    VALERY (Paul)
[1945] *Fonction et mystère de l'Académie »*, *Regards sur le monde actuel et autres essais,* Paris: Gallimard, *1945, p. 250*
    VAPEREAU (Gustave), dir.
[1858] Liouville, *Dictionnaire universel des contemporains*, *I-Z,* Paris, 1858, p.1110.
    VERRIEST (Gustave),
[1934] Évariste Galois et la théorie des équations algébriques, *Revue des questions scientifiques*, mai-juill. 1934.
    VOGT (Paul)
[1982] Identifying Scholarly and Intellectual Communities : A Note on French Philosophy, 1900-1939, *History and Theory*, vol. 21, n° 2, 1982, p. 267-278.
    VOLTERRA (Vito)
[1913] Henri Poincaré: L'oeuvre mathématique, in *Revue du mois*, **15** (1913), 129-154
    WAERDEN (Bartel, van der),
[1972] Die Galoische Theorie von Heinrich Weber bis Emil Artin, *Arch. Hist. Exact. Sci.* 9 (1972), 240-248.
[1985] *A history of Algebra : from Al-Khwàrizmi to Emmy Noether*, New York, Springer Verlag, 1985
    WANTZEL (Pierre Laurent),
[1843] Classification des nombres incommensurables d'origine algébrique, *Nouvelles annales de mathématiques,* 1re série, t. 2 (1843), p. 117-127.

[1845] De l'impossibilité de résoudre toutes les équations algébriques avec des radicaux », *Nouvelles annales de mathématiques*, t. 4, 1845, p.57-66.
  WEBER (Anne-Gaëlle)
[2011] La panthéonisation de la grandeur savante. Les éloges funèbres de l'Académie des sciences de la Belle époque, *in* Weber et al (dir), *Panthéons scientifiques et littéraires*, Presses de l'Université d'Artois, to appear in 2012.
  WEBER (Heinrich),
[1893] Die allgemeinen Grundlagen der Galois'schen Gleichungstheorie, *Mathematische Annalen*, vol. 43, 1893, p. 521-549.
[1895] *Lehrbuch der Algebra*, Braunschweig, F. Vieweg und Sohn, 1895-1896, 2 vols.
  WINTER (Maximilien),
[1908] Importance philosophique de la théorie des nombres, *Revue de métaphysique et de morale,* p. 321-345.
[1910] Caractères de l'algèbre moderne, *Revue de métaphysique et de morale*, année 18, n° 4, 1910, p. 491-529.
[1919] [Recension sur] Pierre Boutroux, *Les principes de l'analyse mathématique*, *Revue de métaphysique et de morale,* 1919, p. 649-667.
  ZERNER (Martin)
[1991] Le règne de Joseph Bertrand (1874-1900), *in* [Gispert 1991, p. 298-322].